\documentclass[12pt]{amsart}
\usepackage{amsxtra}
\usepackage{amssymb, amsmath}
\usepackage{amsthm}
\usepackage {hyperref}

\newtheorem{lemma}{Lemma}
\newtheorem{theorem}{Theorem}
\newtheorem{remark}{Remark}

\newtheorem{property}{Property}
\newtheorem{corollary}{Corollary}

\title[Generalized Divergence Measures]{GENERALIZED SYMMETRIC DIVERGENCE MEASURES AND
INEQUALITIES}
\author{Inder Jeet Taneja}
\address{Inder Jeet Taneja\\
Departamento de Matem\'{a}tica\\
Universidade Federal de Santa Catarina\\
88.040-900 Florian\'{o}polis, SC, Brazil}

\email{taneja@mtm.ufsc.br}

\urladdr{http://www.mtm.ufsc.br/$\sim $taneja}



\begin{document}

\begin{abstract}
In this paper, we have studied the following two
\textit{divergence measures of type s}:
\[
\mathcal V _s (P\vert \vert Q) = \begin{cases}
 {J_s (P\vert \vert Q) = \left[ {s(s - 1)} \right]^{ - 1}\left[
{\sum\limits_{i = 1}^n {\left( {p_i^s q_i^{1 - s} + p_i^{1 - s}
q_i^s }
\right) - 2} } \right],} & {s \ne 0,1} \\
 {J(P\vert \vert Q) = \sum\limits_{i = 1}^n {\left( {p_i - q_i } \right)\ln
\left( {\frac{p_i }{q_i }} \right),} } & {s = 0,1} \\
\end{cases}
\]

\noindent and
\[
\mathcal W _s (P\vert \vert Q) = \begin{cases}
 {IT_s (P\vert \vert Q) = \left[ {s(s - 1)} \right]^{ - 1}\left[
{\sum\limits_{i = 1}^n {\left( {\frac{p_i^{1 - s} + q_i^{1 - s}
}{2}} \right)\left( {\frac{p_i + q_i }{2}} \right)} ^s - 1}
\right],} & {s \ne
0,1} \\
 {I(P\vert \vert Q) = \frac{1}{2}\left[ {\sum\limits_{i = 1}^n {p_i \ln
\left( {\frac{2p_i }{p_i + q_i }} \right) + \sum\limits_{i = 1}^n
{q_i \ln
\left( {\frac{2q_i }{p_i + q_i }} \right)} } } \right],} & {s = 0} \\
 {T(P\vert \vert Q) = \sum\limits_{i = 1}^n {\left( {\frac{p_i + q_i }{2}}
\right)\ln \left( {\frac{p_i + q_i }{2\sqrt {p_i q_i } }} \right)}
,} & {s =
1} \\
\end{cases}
\]

The first measure generalizes the well known \textit{J-divergence}
due to Jeffreys \cite{jef} and Kullback and Leibler \cite{kul}.
The second measure gives a unified generalization of
\textit{Jensen-Shannon divergence} due to Sibson \cite{sib} and
Burbea and Rao \cite{bur1, bur2}, and \textit{arithmetic-geometric
mean divergence }due to Taneja \cite{tan3}. These two measures
contain in particular some well known divergences such as:
\textit{Hellinger's discrimination}, \textit{triangular
discrimination} and \textit{symmetric chi-square divergence}. In
this paper we have studied the properties of the above two
measures and derived some inequalities among them.
\end{abstract}

\maketitle

\section{Introduction}

Let
\[
\Gamma _n = \left\{ {P = (p_1 ,p_2 ,...,p_n )\left| {p_i >
0,\sum\limits_{i = 1}^n {p_i = 1} } \right.} \right\}, \,\, n
\geqslant 2,
\]

\noindent be the set of all complete finite discrete probability
distributions. For all $P,Q \in \Gamma _n $, the following
measures are well known in the literature on information theory
and statistics:\\

\textbf{$\bullet$ Hellinger Discrimination} (Hellinger \cite{hel})
\begin{equation}
\label{eq1} h(P\vert \vert Q) = 1 - B(P\vert \vert Q) =
\frac{1}{2}\sum\limits_{i = 1}^n {(\sqrt {p_i } - \sqrt {q_i }
)^2} \,\, ,
\end{equation}

\noindent where
\begin{equation}
\label{eq2} B(P\vert \vert Q) = \sum\limits_{i = 1}^n \sqrt {p_i
q_i } .
\end{equation}

\noindent is the well-known Bhattacharyya \cite{bha}
\textit{coefficient}.\\

\textbf{$\bullet$ Triangular Discrimination}
\begin{equation}
\label{eq3}
\Delta (P\vert \vert Q) = 2\left[ {1 - W(P\vert \vert Q)} \right] =
\sum\limits_{i = 1}^n {\frac{(p_i - q_i )^2}{p_i + q_i }} ,
\end{equation}

\noindent where
\begin{equation}
\label{eq4} W(P\vert \vert Q) = \sum\limits_{i = 1}^n {\frac{2p_i
q_i }{p_i + q_i }},
\end{equation}

\noindent is the well-known \textit{harmonic mean divergence}.\\

\textbf{$\bullet$ Symmetric Chi-square Divergence} (Dragomir et
al. \cite{dsb})
\begin{equation}
\label{eq5}
\Psi (P\vert \vert Q) = \chi ^2(P\vert \vert Q) + \chi ^2(Q\vert \vert P) =
\sum\limits_{i = 1}^n {\frac{(p_i - q_i )^2(p_i + q_i )}{p_i q_i }} ,
\end{equation}

\noindent where
\begin{equation}
\label{eq6}
\chi ^2(P\vert \vert Q) = \sum\limits_{i = 1}^n {\frac{(p_i - q_i )^2}{q_i
}} = \sum\limits_{i = 1}^n {\frac{p_i^2 }{q_i } - 1} ,
\end{equation}

\noindent is the well-known $\chi ^2 - $\textit{divergence}
(Pearson \cite{pea})\\

\textbf{$\bullet$ J-Divergence} (Jeffreys \cite{jef};
Kullback-Leibler \cite{kul})
\begin{equation}
\label{eq7}
J(P\vert \vert Q) = \sum\limits_{i = 1}^n {(p_i - q_i )\ln (\frac{p_i }{q_i
})} .
\end{equation}

\bigskip
\textbf{$\bullet$ Jensen-Shannon Divergence} (Sibson \cite{sib};
Burbea and Rao \cite{bur1, bur2})
\begin{equation}
\label{eq8}
I(P\vert \vert Q) = \frac{1}{2}\left[ {\sum\limits_{i = 1}^n {p_i \ln \left(
{\frac{2p_i }{p_i + q_i }} \right) + } \sum\limits_{i = 1}^n {q_i \ln \left(
{\frac{2q_i }{p_i + q_i }} \right)} } \right].
\end{equation}

\bigskip
\textbf{$\bullet$ Arithmetic-Geometric Divergence} (Taneja
\cite{tan3})
\begin{equation}
\label{eq9}
T(P\vert \vert Q) = \sum\limits_{i = 1}^n {\left( {\frac{p_i + q_i }{2}}
\right)\ln \left( {\frac{p_i + q_i }{2\sqrt {p_i q_i } }} \right)} .
\end{equation}

After simplification, we can write
\begin{equation}
\label{eq10} J(P\vert \vert Q) = 4\left[ {I(P\vert \vert Q) +
T(P\vert \vert Q)} \right].
\end{equation}

The measures $J(P\vert \vert Q)$, $I(P\vert \vert Q)$ and
$T(P\vert \vert Q)$ can also be written as
\begin{align}
\label{eq11} J(P\vert \vert Q) & = K(P\vert \vert Q) + K(Q\vert
\vert P),\\
\label{eq12} I(P\vert \vert Q) & = \frac{1}{2}\left[ {K\left(
{P\vert \vert \frac{P + Q}{2}} \right) + K\left( {Q\vert \vert
\frac{P + Q}{2}} \right)} \right]\\
\intertext{and} \label{eq13} T(P\vert \vert Q) & =
\frac{1}{2}\left[ {K\left( {\frac{P + Q}{2}\vert \vert P} \right)
+ K\left( {\frac{P + Q}{2}\vert \vert Q} \right)} \right],
\end{align}

\noindent where
\begin{equation}
\label{eq14}
K(P\vert \vert Q) = \sum\limits_{i = 1}^n {p_i \log \left( {\frac{p_i }{q_i
}} \right)} ,
\end{equation}

\noindent is the well known Kullback-Leibler \cite{kul}
\textit{relative information}.

The measure (\ref{eq9}) is also known by \textit{Jensen difference
divergence measure} (Burbea and Rao \cite{bur1, bur2}). The
measure (\ref{eq10}) is new in the literature and is studied for
the first time by Taneja \cite{tan3} and is called
\textit{arithmetic and geometric mean divergence measure.} For
simplicity, the three measures appearing in (\ref{eq8}),
(\ref{eq9}) and (\ref{eq10}) we shall call,
\textit{JS-divergence}, \textit{J-divergence} and the
\textit{AG-divergence} respectively. More details on these
divergence measures can be seen in on line book by Taneja
\cite{tan4}.

We call the measures given in (\ref{eq1}), (\ref{eq3}),
(\ref{eq5}), (\ref{eq7}), (\ref{eq9}) and (\ref{eq10}) by
\textit{symmetric divergence measures}, since they are symmetric
with respect to the probability distributions $P$ and $Q$. While
the measures (\ref{eq6}) and (\ref{eq14}) are not symmetric with
respect to probability distributions.

\section{Generalizations of Symmetric Divergence
Measures}

In this section, we shall present new generalizations of the
symmetric divergence measures given in Section 1. Before that,
first we shall present a well known generalization of
Kullback-Leibler's \textit{relative information}.\\

\textbf{$\bullet$ Relative Information of Type s}
\begin{equation}
\label{eq15} \Phi _s (P\vert \vert Q) = \begin{cases}
 {K_s (P\vert \vert Q) = \left[ {s(s - 1)} \right]^{ - 1}\left[
{\sum\limits_{i = 1}^n {p_i^s q_i^{1 - s} - 1} } \right],} & {s \ne 0,1}
\\\\
 {K(Q\vert \vert P) = \sum\limits_{i = 1}^n {q_i \ln \left( {\frac{q_i }{p_i
}} \right)} ,} & {s = 0} \\\\
 {K(P\vert \vert Q) = \sum\limits_{i = 1}^n {p_i \ln \left( {\frac{p_i }{q_i
}} \right)} ,} & {s = 1} \\
\end{cases},
\end{equation}

\noindent for all $s \in \mathbb{R}$.\\

The measure (\ref{eq15}) is due to Cressie and Read \cite{crr}.
For more studies on this measure refer to Taneja \cite{tan5} and
Taneja and Kumar \cite{tak, kut} and reference therein.

The measure \ref{eq15} admits the following particular cases:
\begin{itemize}
\item[(i)] $ \Phi _{ - 1} (P\vert \vert Q) = \frac{1}{2}\chi
^2(Q\vert \vert P).$

\item[(ii)] $\Phi _0 (P\vert \vert Q) = K(Q\vert \vert P).$

\item[(iii)] $\Phi _{1 / 2} (P\vert \vert Q) = 4\left[ {1 -
B(P\vert \vert Q)} \right] = 4h(P\vert \vert Q).$

\item[(iv)] $\Phi _1 (P\vert \vert Q) = K(P\vert \vert Q).$

\item[(v)] $\Phi _2 (P\vert \vert Q) = \frac{1}{2}\chi ^2 (P\vert
\vert Q).$
\end{itemize}

Here, we observe that $\Phi _2 (P\vert \vert Q) = \Phi _{ - 1}
(Q\vert \vert P)$ and $\Phi _1 (P\vert \vert Q) = \Phi _0 (Q\vert
\vert P)$.\\

\subsection*{2.1. J-Divergence of Type s}

Replace $K(P\vert \vert Q)$ by $\Phi _s (P\vert \vert Q)$ in the
relation (\ref{eq11}), we get
\begin{align}
\label{eq16} \mathcal V _s (P\vert \vert Q) & = \Phi _s (P\vert
\vert
Q) + \Phi _s (Q\vert \vert P)\\
& = \begin{cases}
 {J_s (P\vert \vert Q) = \left[ {s(s - 1)} \right]^{ - 1}\left[
{\sum\limits_{i = 1}^n {\left( {p_i^s q_i^{1 - s} + p_i^{1 - s} q_i^s }
\right) - 2} } \right],} & {s \ne 0,1} \\\\
 {J(P\vert \vert Q) = \sum\limits_{i = 1}^n {\left( {p_i - q_i } \right)\ln
\left( {\frac{p_i }{q_i }} \right),} } & {s = 0,1} \\
\end{cases}.\notag
\end{align}

The expression (\ref{eq16}) admits the following particular cases:
\begin{itemize}
\item[(i)] $\mathcal V _{ - 1} (P\vert \vert Q) = \mathcal V _2
(P\vert \vert Q) = \frac{1}{2}\Psi (P\vert \vert Q).$

\item[(ii] $\mathcal V _0 (P\vert \vert Q) = \mathcal V _1 (P\vert
\vert Q) = J(P\vert \vert Q).$

\item[(iii)] $\mathcal V _{1 / 2} (P\vert \vert Q) = 8\,h(P\vert
\vert Q).$
\end{itemize}

\begin{remark} The expression (\ref{eq16}) is the modified form of the
measure already known in the literature:
\begin{equation}
\label{eq175} \mathcal V _s^1 (P\vert \vert Q) = \begin{cases}
 {J_s (P\vert \vert Q) = (s - 1)^{ - 1}\left[
{\sum\limits_{i = 1}^n {\left( {p_i^s q_i^{1 - s} + p_i^{1 - s}
q_i^s }
\right) - 2} } \right],} & {s \ne 1}, s > 0 \\\\
 {J(P\vert \vert Q) = \sum\limits_{i = 1}^n {\left( {p_i - q_i } \right)\ln
\left( {\frac{p_i }{q_i }} \right),} } & {s = 1} \\
\end{cases}.
\end{equation}

For the propertied of the measure (\ref{eq175}) refer to Burbea
and Rao \cite{bur1, bur2}, Taneja \cite{tan2, tan3, tan4}, etc.
For the axiomatic characterization of this measure refer to Rathie
and Sheng \cite{ras} and Taneja \cite{tan1}. The measures
(\ref{eq16}) considered here differs in constant and it permits in
considering negative values of the parameter $s$.
\end{remark}

\subsection*{2.2. Unified AG and JS -- Divergence of Type
s}

Replace $K(P\vert \vert Q)$ by $\Phi _s (P\vert \vert Q)$ in the
relation (\ref{eq13}) interestingly we have a unified
generalization of the \textit{AG and JS -- divergence} given by

\begin{align}
\label{eq17} \mathcal W _s (P\vert \vert Q) & = \frac{1}{2}\left[
{\Phi _s \left( {\frac{P + Q}{2}\vert \vert P} \right) + \Phi _s
\left(
{\frac{P + Q}{2}\vert \vert Q} \right)} \right]\\
& = \begin{cases}
 {IT_s (P\vert \vert Q) = \left[ {s(s - 1)} \right]^{ - 1}\left[
{\sum\limits_{i = 1}^n {\left( {\frac{p_i^{1 - s} + q_i^{1 - s} }{2}}
\right)\left( {\frac{p_i + q_i }{2}} \right)} ^s - 1} \right],} & {s \ne
0,1} \\\\
 {I(P\vert \vert Q) = \frac{1}{2}\left[ {\sum\limits_{i = 1}^n {p_i \ln
\left( {\frac{2p_i }{p_i + q_i }} \right) + \sum\limits_{i = 1}^n {q_i \ln
\left( {\frac{2q_i }{p_i + q_i }} \right)} } } \right],} & {s = 0}
\\\\
 {T(P\vert \vert Q) = \sum\limits_{i = 1}^n {\left( {\frac{p_i + q_i }{2}}
\right)\ln \left( {\frac{p_i + q_i }{2\sqrt {p_i q_i } }} \right)} ,} & {s =
1} \\
\end{cases}.\notag
\end{align}

The measure (\ref{eq17}) admits the following particular cases:
\begin{itemize}
\item[(i)] $\mathcal W _{ - 1} (P\vert \vert Q) =
\frac{1}{4}\Delta (P\vert \vert Q)$.

\item[(ii)] $\mathcal W _0 (P\vert \vert Q) = I(P\vert \vert Q)$.

\item[(iii)] $\mathcal W _{1/2} (P\vert \vert Q) = 4\,d(P\vert
\vert Q)$.

\item[(iv)] $\mathcal W _1 (P\vert \vert Q) = T(P\vert \vert Q)$.

\item[(v)] $\mathcal W _2 (P\vert \vert Q) = \frac{1}{16}\Psi
(P\vert \vert Q)$.
\end{itemize}

The measure $d(P\vert \vert Q)$ given in part (iii) is not studied
elsewhere and given by
\begin{equation}
\label{eq181} d(P\vert \vert Q) = 1 - \sum\limits_{i = 1}^n
{\left( {\frac{\sqrt {p_i } + \sqrt {q_i } }{2}} \right)} \left(
{\sqrt {\frac{p_i + q_i }{2}} } \right).
\end{equation}

A relation of the measure (\ref{eq181}) with \textit{Hellinger's
discrimination} is given in the last section. Connections of the
measure (\ref{eq181}) with mean divergence measures can be seen in
Taneja \cite{tan8}.

\bigskip
 We can also write
\begin{equation}
\label{eq18} \mathcal W _{1 - s} (P\vert \vert Q) =
\frac{1}{2}\left[ {\Phi _s \left( {P\vert \vert \frac{P + Q}{2}}
\right) + \Phi _s \left( {Q\vert \vert \frac{P + Q}{2}} \right)}
\right].
\end{equation}

Thus we have two \textit{symmetric divergences of type s} given by
(\ref{eq16}) and (\ref{eq17}) generalizing the six
\textit{symmetric divergence measures} given in Section 1. In this
paper our aim is to study the \textit{symmetric divergences of
type s} and to find inequalities among them. These studies we
shall do by making use of the properties of Csisz\'{a}r's
\textit{f-divergence}.

\section{Csisz\'{a}r's $f-$Divergence and Its
Properties}

Given a function $f:(0,\infty ) \to \mathbb{R}$, the
\textit{f-divergence} measure introduced by Csisz\'{a}r's
\cite{csi1} is given by
\begin{equation}
\label{eq19}
C_f (P\vert \vert Q) =
\sum\limits_{i = 1}^n {q_i f\left( {\frac{p_i }{q_i }} \right)} ,
\end{equation}

\noindent for all $P,Q \in \Gamma _n $.

The following theorem is well known in the literature \cite{csi1,
csi2}.

\begin{theorem} \label{the31} If the function
$f$ is convex and normalized, i.e., $f(1) = 0$, then the
\textit{f-divergence}, $C_f (P\vert \vert Q)$ is
\textit{nonnegative} and \textit{convex} in the pair of
probability distribution $(P,Q) \in \Gamma _n \times \Gamma _n $.
\end{theorem}

The theorem given below give bounds on the measure (\ref{eq19}).

\begin{theorem} \label{the32} Let $f:\mathbb{R}_ + \to \mathbb{R}$ be
differentiable convex and normalized i.e., $f(1) = 0$. If $P,Q \in
\Gamma _n $, are such that $0 < r \leqslant \frac{p_i }{q_i }
\leqslant R < \infty $, $\forall i \in \{1,2,...,n\}$, for some
$r$ and $R$ with $0 < r \leqslant 1 \leqslant R < \infty $, $r \ne
R$, then we have
\begin{equation}
\label{eq20}
0 \leqslant C_f (P\vert \vert Q) \leqslant E_{C_f } (P\vert \vert Q)
\leqslant A_{C_f } (r,R),
\end{equation}

\noindent and
\begin{equation}
\label{eq21}
0 \leqslant C_f (P\vert \vert Q) \leqslant B_{C_f } (r,R) \leqslant A_{C_f }
(r,R),
\end{equation}

\noindent where
\begin{equation}
\label{eq22}
E_{C_f } (P\vert \vert Q) = \sum\limits_{i = 1}^n {(p_i - q_i )}
{f}'(\frac{p_i }{q_i }),
\end{equation}
\begin{equation}
\label{eq23} A_{C_f } (r,R) = \frac{1}{4}(R - r)\left(
{\,\,{f}'(R) - {f}'(r)} \right)
\end{equation}

\noindent and
\begin{equation}
\label{eq24}
B_{C_f } (r,R) = \frac{(R - 1)f(r) + (1 - r)f(R)}{R - r}.
\end{equation}
\end{theorem}

The proof is based on the following lemma due to Dragomir
\cite{dra1}.

\begin{lemma} \label{lem31} Let $f:I \subset \mathbb{R}_ + \to
\mathbb{R}$ be a differentiable convex function on the interval
$I$, $x_i \in \mathop I\limits^o $ ($\mathop I\limits^o $ is the
interior of $I)$, $\lambda _i \geqslant 0$ ($i = 1,2,...,n)$ with
$\sum\limits_{i = 1}^n {\lambda _i = 1} $. If $m,\,\,M \in \mathop
I\limits^o $ and $m \leqslant x_i \leqslant M$, $\forall i =
1,2,...,n$, then we have the inequalities:
\begin{align}
\label{eq25} 0 & \leqslant \sum\limits_{i = 1}^n {\lambda _i f(x_i
) - f\left( {\sum\limits_{i = 1}^n {\lambda _i x_i } } \right)}\\
& \leqslant \sum\limits_{{\i} = 1}^n {\lambda _i x_i {f}'(x_i ) -
\left( {\sum\limits_{i = 1}^n {\lambda _i x_i } } \right)\left(
{\sum\limits_{i = 1}^n {\lambda _i {f}'(x_i )} } \right)}\notag\\
& \leqslant \frac{1}{4}(M - m)\left( {{f}'(M) - {f}'(m)}
\right).\notag
\end{align}
\end{lemma}

As a consequence of above theorem we have the following corollary.

\begin{corollary} \label{cor31} For all $a,b,\upsilon ,\omega \in
(0,\infty )$, the following inequalities hold:
\begin{align}
\label{eq26} 0 & \leqslant \frac{\upsilon f(a) + \omega
f(b)}{\upsilon + \omega } - f\left( {\frac{\upsilon a + \omega
b}{\upsilon + \omega }} \right)\\
& \leqslant \frac{\upsilon a{f}'(a) + \omega b{f}'(b)}{\upsilon +
\omega } - \left( {\frac{\upsilon a + \omega b}{\upsilon + \omega
}} \right)\left( {\frac{\upsilon {f}'(a) + \omega
{f}'(b)}{\upsilon + \omega }} \right)\notag\\
& \leqslant \frac{1}{4}(b - a)\left( {{f}'(b) - {f}'(a)}
\right).\notag
\end{align}
\end{corollary}

\begin{proof} It follows from Lemma \ref{lem31}, by taking $\lambda
_1 = \frac{\upsilon }{\upsilon + \omega }$, $\lambda _2 =
\frac{\omega }{\upsilon + \omega }$, $\lambda _3 = ... = \lambda
_n = 0$, $x_1 = a$, $x_2 = b$, $x_2 = ... = x_n = 0$.
\end{proof}

\begin{proof} \textit{of the Theorem \ref{the32}}. For all $P,Q \in \Gamma _n
$, take $x = \frac{p_i }{q_i }$ in (\ref{eq25}), $\lambda _i = q_i
$ and sum over all $i = 1,2,...,n$ we get the inequalities
(\ref{eq20}).

Again, take $\upsilon = R - x$, $\omega = x - r$, $a = r$ and $b =
R$ in (\ref{eq26}), we get
\begin{align}
\label{eq27} 0 & \leqslant \frac{(R - x){f}(r) + (x -
r){f}(R)}{R - r} - f(x)\\
& \leqslant \frac{(R - x)(x - r)}{R - r}\left[ {{f}'(R) - {f}'(r)}
\right]\notag \\
& \leqslant \frac{1}{4}(R - r)\left( {{f}'(R) - {f}'(r)}
\right).\notag
\end{align}

From the first part of the inequalities (\ref{eq27}), we get
\begin{equation}
\label{eq28} f(x) \leqslant \frac{(R - x){f}(r) + (x - r){f}(R)}{R
- r}.
\end{equation}

For all $P,Q \in \Gamma _n $, take $x = \frac{p_i }{q_i }$ in
(\ref{eq27}) and (\ref{eq28}), multiply by $q_i $ and sum over all
$i = 1,2,...,n$, we get
\begin{align}
\label{eq29} 0 & \leqslant B_{C_f } (r,R) - C_f (P\vert \vert Q)
\\
& \leqslant \frac{(R - 1)(1 - r)\left( {{f}'(R) - {f}'(r)}
\right)}{R - r} \leqslant A_{C_f } (r,R)\notag
\end{align}

\noindent and
\begin{equation}
\label{eq30} 0 \leqslant C_f (P\vert \vert Q) \leqslant B_{C_f }
(r,R),
\end{equation}

\noindent respectively.

The expression (\ref{eq30}) completes the $l.h.s.$ of the
inequalities (\ref{eq21}). In order to prove $r.h.s.$ of the
inequalities (\ref{eq21}), let us take $x = 1$ in (\ref{eq27}) and
use the fact that $f(1) = 0$, we get
\begin{equation}
\label{eq31} 0 \leqslant B_{C_f } (r,R)
 \leqslant \frac{(R - 1)(1 - r)\left( {{f}'(R) - {f}'(r)} \right)}{R - r}
 \leqslant A_{C_f } (r,R).
\end{equation}

From (\ref{eq31}), we conclude the $r.h.s.$ of the inequalities
(\ref{eq21}).
\end{proof}

\begin{remark}
We observe that the inequalities (\ref{eq20}) and (\ref{eq21}) are
the improvement over Dragomir's \cite{dra2, dra3} work. From the
inequalities (\ref{eq21}) and (\ref{eq31}) we observe that there
is better bound for $B_{C_f } (r,R)$ instead of $A_{C_f } (r,R)$.
From the $r.h.s.$ of the inequalities (\ref{eq31}), we conclude
the following inequality among $r$ and $R$:
\begin{equation}
\label{eq32} (R - 1)(1 - r) \leqslant \frac{1}{4}(R - r)^2.
\end{equation}
\end{remark}

\begin{theorem} \label{the33} $($Dragomir et al. \cite{dgp1, dgp2}$)$.
$(i)$ Let $P,Q \in \Gamma _n $ be such that $0 < r \leqslant
\frac{p_i }{q_i } \leqslant R < \infty $, $\forall i \in
\{1,2,...,n\}$, for some $r$ and $R$ with $0 < r \leqslant 1
\leqslant R < \infty $, $r \ne R$. Let $f:[0,\infty ) \to
\mathbb{R}$ be a normalized mapping, i.e., $f(1) = 0$ such that
${f}'$ is locally absolutely continuous on $[r,R]$ and there
exists $\alpha ,\,\,\beta $ satisfying
\begin{equation}
\label{eq33} \alpha \leqslant {f}''(x) \leqslant \beta , \,\,
\forall x \in (r,R).
\end{equation}

Then
\begin{equation}
\label{eq34} \left| {C_f (P\vert \vert Q) - \frac{1}{2}E_{C_f }
(P\vert \vert Q)} \right| \leqslant \frac{1}{8}(\beta - \alpha
)\chi ^2(P\vert \vert Q)
\end{equation}

\noindent and
\begin{equation}
\label{eq35} \left| {C_f (P\vert \vert Q) - E_{C_f }^\ast (P\vert
\vert Q)} \right| \leqslant \frac{1}{8}(\beta - \alpha )\chi
^2(P\vert \vert Q),
\end{equation}

\noindent where $E_{C_f } (P\vert \vert Q)$ is as given by
(\ref{eq22}), $\chi ^2(P\vert \vert Q)$ is as given by (\ref{eq6})
and
\begin{equation}
\label{eq36} E_{C_f }^\ast (P\vert \vert Q) = 2\,E_{C_f }
\left(\frac{P+Q}{2}\vert \vert Q\right) = \sum\limits_{i = 1}^n
{(p_i - q_i ){f}'\left( {\frac{p_i + q_i }{2q_i }} \right)} .
\end{equation}

$(ii)$ Additionally, if $f:[r,R] \to \mathbb{R}$ with ${f}'''$
absolutely continuous on $[r,R]$ and ${f}''' \in L_\infty [r,R]$,
then
\begin{equation}
\label{eq37} \left| {C_f (P\vert \vert Q) - \frac{1}{2}E_{C_f }
(P\vert \vert Q)} \right| \leqslant \frac{1}{12}\left\| {f}'''
\right\|_\infty \vert \chi \vert ^3(P\vert \vert Q)
\end{equation}

\noindent and
\begin{equation}
\label{eq38} \left| {C_f (P\vert \vert Q) - E_{C_f }^\ast (P\vert
\vert Q)} \right| \leqslant \frac{1}{24}\left\| {f}'''
\right\|_\infty \vert \chi \vert ^3(P\vert \vert Q),
\end{equation}

\noindent where
\begin{equation}
\label{eq39} \vert \chi \vert ^3(P\vert \vert Q) = \sum\limits_{i
= 1}^n {\frac{\vert p_i - q_i \vert ^3}{q_i^2 }}
\end{equation}

\noindent and
\begin{equation}
\label{eq40} \left\| {f}''' \right\|_\infty = ess\mathop {\sup
}\limits_{x \in [r,R]} \vert {f}'''\vert .
\end{equation}
\end{theorem}

\begin{theorem} \label{the34} $($Dragomir et al. \cite{dgp3}$)$. Suppose
$f:[r,R] \to \mathbb{R}$ is differentiable and ${f}'$ is of
bounded variation, i.e., $\mathop V\limits_r^R ({f}') = \int_r^R
{\vert {f}''(t)\vert dt < \infty } $. Let the constants $r,R$
satisfy the conditions:
\begin{itemize}
\item[(i)] $0 < r < 1 < R < \infty $;
\item[(ii)] $0 < r \leqslant
\frac{p_i }{q_i } \leqslant R < \infty $, for $i = 1,2,...,n$.
\end{itemize}
Then
\begin{equation}
\label{eq41} \left| {C_f (P\vert \vert Q) - \frac{1}{2}E_{C_f }
(P\vert \vert Q)} \right| \leqslant \mathop V\limits_r^R
({f}')V(P\vert \vert Q)
\end{equation}

\noindent and
\begin{equation}
\label{eq42} \left| {C_f (P\vert \vert Q) - E_{C_f }^\ast (P\vert
\vert Q)} \right| \leqslant \frac{1}{2}\mathop V\limits_r^R
({f}')V(P\vert \vert Q),
\end{equation}

\noindent where
\begin{equation}
\label{eq43} V(P\vert \vert Q) = \sum\limits_{i = 1}^n {\left|
{p_i - q_i } \right|} .
\end{equation}
\end{theorem}

\begin{remark} $(i)$ If the third order derivative of $f$ exists and let us
suppose that it is either positive or negative. Then the function
${f}''$ is either monotonically increasing or decreasing. In view
of this we can write
\begin{equation}
\label{eq44} \beta - \alpha = k(f)\left[ {{f}''(R) - {f}''(r)}
\right],
\end{equation}
\noindent where
\begin{equation}
\label{eq45} k(f) = \begin{cases}
 { - 1,} & {\mbox{if }{f}''\mbox{ is monotonically decresing}} \\
 {1,} & {\mbox{if }{f}''\mbox{ is monotonically increasing}} \\
\end{cases}.
\end{equation}

$(ii)$ Let the function $f(x)$ considered in the Theorem 3.4 be
convex in $(0,\infty )$, then ${f}''(x) \geqslant 0$. This gives
\begin{align}
\label{eq46} \mathop V\limits_r^R ({f}') & = \int_r^R {\vert
{f}''(t)\vert dt} \\
& = \int_r^R {{f}''(t)dt = {f}'(R) - {f}'(r)}\notag\\
& = \frac{4}{R - r}A_{C_f } (r,R).\notag
\end{align}

Under these considerations, the bounds (\ref{eq41}) and
(\ref{eq42}) can be re-written as
\begin{equation}
\label{eq47} \left| {C_f (P\vert \vert Q) - \frac{1}{2}E_{C_f }
(P\vert \vert Q)} \right| \leqslant \frac{4}{R - r}A_{C_f }
(r,R)V(P\vert \vert Q)
\end{equation}
\noindent and
\begin{equation}
\label{eq48} \left| {C_f (P\vert \vert Q) - E_{C_f }^\ast (P\vert
\vert Q)} \right| \leqslant \frac{2}{R - r}A_{C_f } (r,R)V(P\vert
\vert Q).
\end{equation}
\end{remark}

Based on above remarks we can restate and combine the Theorems
\ref{the33} and \ref{the34}.

\begin{theorem} \label{the35} Let $P,Q \in \Gamma _n $ be such that $0
< r \leqslant \frac{p_i }{q_i } \leqslant R < \infty $, $\forall i
\in \{1,2,...,n\}$, for some $r$ and $R$ with $0 < r < 1 < R <
\infty $. Let $f:\mathbb{R}_ + \to \mathbb{R}$ be differentiable
convex, normalized, of bounded variation, and second derivative is
monotonic with ${f}'''$ absolutely continuous on $[r,R]$ and
${f}''' \in L_\infty [r,R]$, then

\begin{align}
\label{eq49} & \left| {C_f (P\vert \vert Q) - \frac{1}{2}E_{C_f }
(P\vert \vert Q)} \right|\\
& \,\, \leqslant \min \left\{ {\frac{1}{8}k(f)\left[ {{f}''(R) -
{f}''(r)} \right]\chi ^2(P\vert \vert Q),} \right.\notag\\
& \qquad \qquad \left. {\frac{1}{12}\left\| {f}''' \right\|_\infty
\vert \chi \vert ^3(P\vert \vert Q),\mbox{ }\left[ {{f}'(R) -
{f}'(r)} \right]V(P\vert \vert Q)} \right\},\notag
\end{align}

\noindent and
\begin{align}
\label{eq50} & \left| {C_f (P\vert \vert Q) - E_{C_f }^\ast
(P\vert \vert Q)} \right|\\
& \,\, \leqslant \min \left\{ {\frac{1}{8}k(f)\left[ {{f}''(R) -
{f}''(r)} \right]\chi ^2(P\vert \vert Q),} \right.\notag\\
& \qquad \qquad \left. {\frac{1}{24}\left\| {f}''' \right\|_\infty
\vert \chi \vert ^3(P\vert \vert Q),\mbox{ }\frac{1}{2}\left[
{{f}'(R) - {f}'(r)} \right]V(P\vert \vert Q)} \right\},\notag
\end{align}

\noindent where $k(f)$ is as given by (\ref{eq45}).
\end{theorem}

\begin{remark}
The measures (\ref{eq6}), (\ref{eq39}) and (\ref{eq43}) are the
particular cases of Vajda \cite{vaj} $\vert \chi \vert ^m -
$\textit{divergence} given by
\begin{equation}
\label{eq51} \vert \chi \vert ^m(P\vert \vert Q) = \sum\limits_{i
= 1}^n {\frac{\vert p_i - q_i \vert ^m}{q_i^{m - 1} }} ,\mbox{ }m
\geqslant 1.
\end{equation}

The above measure (\ref{eq51}) \cite{ced} \cite{dgp1} satisfies
the following properties:
\begin{align}
\label{eq52} \vert \chi \vert ^m(P\vert \vert Q) & \leqslant
\frac{(1 - r)(R - 1)}{(R - r)}\left[ {(1 - r)^{m - 1} + (R - 1)^{m
- 1}} \right] \\
& \leqslant \left( {\frac{R - r}{2}} \right)^m, \,\, m \geqslant
1\notag
\end{align}

\noindent and
\begin{equation}
\label{eq53} \left( {\frac{1 - r^m}{1 - r}} \right)V(P\vert \vert
Q) \leqslant \vert \chi \vert ^m(P\vert \vert Q) \leqslant \left(
{\frac{R^m - 1}{R - 1}} \right)V(P\vert \vert Q), \,\, m \geqslant
1.
\end{equation}

Take $m=2, 3$ and $1$, in (\ref{eq52}), we get
\begin{equation}
\label{eq54} \chi ^2(P\vert \vert Q) \leqslant (R - 1)(1 - r)
\leqslant \frac{(R - r)^2}{4},
\end{equation}
\begin{equation}
\label{eq55} \vert \chi \vert ^3(P\vert \vert Q) \leqslant
\frac{1}{2}\frac{(R - 1)(1 - r)}{R - r}\left[ {(1 - r)^2 + (R -
1)^2} \right] \leqslant \frac{1}{8}(R - r)^3
\end{equation}

\noindent and
\begin{equation}
\label{eq56} V(P\vert \vert Q) \leqslant \frac{2(R - 1)(1 - r)}{(R
- r)} \leqslant \frac{1}{2}(R - r).
\end{equation}
\noindent respectively.

In view of the last inequalities given in (\ref{eq54}),
(\ref{eq55}) and (\ref{eq56}), the bounds given in (\ref{eq49})
and (\ref{eq50}) can be written in terms of $r,R$ as

\begin{align}
\label{eq57} & \left| {C_f (P\vert \vert Q) - \frac{1}{2}E_{C_f }
(P\vert \vert
Q)} \right|\\
& \,\, \leqslant \frac{(R - r)^2}{4}\min \left\{
{\frac{1}{8}k(f)\left[
{{f}''(R) - {f}''(r)} \right],} \right. \notag\\
& \qquad \qquad  \qquad  \qquad  \left. {\frac{R - r}{24}\left\|
{f}''' \right\|_\infty ,\,\,\frac{\mbox{2}\left[ {{f}'(R) -
{f}'(r)} \right]}{R - r}} \right\}\notag
\end{align}

\noindent and
\begin{align}
\label{eq58} & \left| {C_f (P\vert \vert Q) - E_{C_f }^\ast
(P\vert \vert Q)} \right|\\
&\,\, \leqslant \frac{(R - r)^2}{4}\min \left\{
{\frac{1}{8}k(f)\left[ {{f}''(R) - {f}''(r)} \right],}
\right.\notag\\
& \qquad \qquad \qquad \qquad \left. {\frac{R - r}{48}\left\|
{f}''' \right\|_\infty ,\,\,\frac{{f}'(R) - {f}'(r)}{R - r}}
\right\} ,
\end{align}

\noindent respectively.

We observe that the bounds (\ref{eq57}) and (\ref{eq58}) are based
on the \textit{first}, \textit{second} and \textit{third}
\textit{order derivatives} of the \textit{generating function}.
\end{remark}

\begin{theorem} \label{the36} Let $f_1 ,f_2 :I \subset \mathbb{R}_ +
\to \mathbb{R}$ two generating mappings are normalized, i.e., $f_1
(1) = f_2 (1) = 0$ and satisfy the assumptions:
\begin{itemize}
\item[(i)] $f_1 $ and $f_2 $ are twice differentiable on $(r,R)$;

\item[(ii)] there exists the real constants $m,M$ such that $m <
M$ and
\end{itemize}
\begin{equation}
\label{eq59} m \leqslant \frac{f_1 ^{\prime \prime }(x)}{f_2
^{\prime \prime }(x)} \leqslant M, \,\, f_2 ^{\prime \prime }(x)
> 0, \,\, \forall x \in (r,R),
\end{equation}

\noindent then we have
\begin{equation}
\label{eq60} m\,\,C_{f_2 } (P\vert \vert Q) \leqslant C_{f_1 }
(P\vert \vert Q) \leqslant M\,\,C_{f_2 } (P\vert \vert Q)
\end{equation}
\end{theorem}

\begin{proof} Let us consider two functions
\begin{equation}
\label{eq61} \eta _m (x) = f_1 (x) - m\,\,f_2 (x),
\end{equation}

\noindent and
\begin{equation}
\label{eq62} \eta _M (x) = M\,\,f_2 (x) - f_1 (x),
\end{equation}

\noindent where $m$ and $M$ are as given by (\ref{eq59})

Since $f_1 (1) = f_2 (1) = 0$, then $\eta _m (1) = \eta _M (1) =
0$. Also, the functions $f_1 (x)$ and $f_2 (x)$ are twice
differentiable. Then in view of (\ref{eq59}), we have
\begin{equation}
\label{eq63} {\eta }''_m (x) = f_1 ^{\prime \prime }(x) - m\,\,f_2
^{\prime \prime }(x)
 = f_2 ^{\prime \prime }(x)\left( {\frac{f_1 ^{\prime \prime }(x)}{f_2
^{\prime \prime }(x)} - m} \right) \geqslant 0,
\end{equation}

\noindent and
\begin{equation}
\label{eq64} {\eta }''_M (x) = M\,\,f_2 ^{\prime \prime }(x) - f_1
^{\prime \prime }(x)
 = f_2 ^{\prime \prime }(x)\left( {M - \frac{f_1 ^{\prime \prime }(x)}{f_2
^{\prime \prime }(x)}} \right) \geqslant 0,
\end{equation}

\noindent for all $x \in (r,R)$.

In view of (\ref{eq63}) and (\ref{eq64}), we can say that the
functions $\eta _m (x)$ and $\eta _M (x)$ are convex on $(r,R)$.

According to Theorem \ref{the31}, we have
\begin{equation}
\label{eq65} C_{\eta _m } (P\vert \vert Q) = C_{f_1 - mf_2 }
(P\vert \vert Q) = C_{f_1 } (P\vert \vert Q) - m\,\,C_{f_2 }
(P\vert \vert Q) \geqslant 0,
\end{equation}

\noindent and
\begin{equation}
\label{eq66} C_{\eta _M } (P\vert \vert Q) = C_{Mf_2 - f_1 }
(P\vert \vert Q) = M\,\,C_{f_2 } (P\vert \vert Q) - C_{f_1 }
(P\vert \vert Q) \geqslant 0.
\end{equation}

Combining (\ref{eq65}) and (\ref{eq66}) we get (\ref{eq60}).
\end{proof}

For futher properties of the measure (\ref{eq19}) based on the
conditions of Theorem \ref{the36} refer to Taneja \cite{tan8}.

\begin{remark} $(i)$ From now onwards, unless otherwise specified,
it is understood that, if there are $r,R$, then $0 < r \leqslant
\frac{p_i }{q_i } \leqslant R < \infty $, $\forall i \in
\{1,2,...,n\}$, with $0 < r < 1 < R < \infty $, $P = (p_1 ,p_2
,....,p_n ) \in \Gamma _n $ and $Q = (q_1 ,q_2 ,....,q_n ) \in
\Gamma _n $.

$(ii)$ In some particular cases studied below, we shall use the
\textit{p-logarithmic power mean} \cite{sto} given by
\begin{equation}
\label{eq67} L_p (a,b) =\begin{cases}
 {\left[ {\frac{b^{p + 1} - a^{p + 1}}{(p + 1)(b - a)}}
\right]^{\frac{1}{p}},} & {p \ne - 1,0} \\\\
 {\frac{b - a}{\ln b - \ln a},} & {p = - 1} \\\\
 {\frac{1}{e}\left[ {\frac{b^b}{a^a}} \right]^{\frac{1}{b - a}},} & {p = 0}
\\
\end{cases},
\end{equation}

\noindent for all $p \in \mathbb{R}$, $a \ne b$. In particular,
we shall use the following notation
\begin{equation}
\label{eq68} L_p^p (a,b) =\begin{cases}
 {\frac{b^{p + 1} - a^{p + 1}}{(p + 1)(b - a)},} & {p \ne - 1,0}
 \\\\
 {\frac{\ln b - \ln a}{b - a},} & {p = - 1} \\\\
 {1,} & {p = 0} \\
\end{cases},
\end{equation}

\noindent for all $p \in \mathbb{R}$, $a \ne b$.
\end{remark}

\section{Bounds on Generalized Divergence Measures}

In this section we shall show that the generalized measures given
in Section 2 are the particular cases of the Csisz\'{a}r's
\textit{f-divergence}. Also, we shall give bounds on these
measures similar to Theorems \ref{the31}-\ref{the35}. The
applications of Theorem \ref{the36} are given in Section 5.

\subsection*{4.1. Bounds on J-Divergence of Type \textit{s}}
Initially we shall give two important properties of
\textit{J-divergence of type s}.

\begin{property} \label{prt41} The measure $\mathcal V _s (P\vert \vert Q)$
is \textit{nonnegative} and \textit{convex} in the pair of
probability distributions $(P,Q) \in \Gamma _n \times \Gamma _n $
for all $s \in ( - \infty ,\infty )$.
\end{property}

\begin{proof} For all $x > 0$ and $s \in ( - \infty ,\infty
)$, let us consider in (\ref{eq15}),
\begin{equation}
\label{eq69} \phi _s (x) =\begin{cases}
 {\left[ {s(s - 1)} \right]^{ - 1}\left[ {x^s + x^{1 - s} - (1 + x)}
\right],} & {s \ne 0,1} \\
 {(x - 1)\ln x,} & {s = 0,1} \\
\end{cases},
\end{equation}

\noindent then we have $C_f (P\vert \vert Q)=\mathcal V _s \left(
{P\vert \vert Q} \right)$, where $\mathcal V _s \left( {P\vert
\vert Q} \right)$ is given by (\ref{eq16}).

Moreover,
\begin{equation}
\label{eq70} \phi _s ^\prime (x) =\begin{cases}
 {\left[ {s(s - 1)} \right]^{ - 1}\left[ {s(x^{s - 1} + x^{ - s}) + x^{ - s}
- 1} \right],} & {s \ne 0,1} \\
 {1 - x^{ - 1} + \ln x,} & {s = 0,1} \\
\end{cases},
\end{equation}

\noindent and
\begin{equation}
\label{eq71} \phi _s ^{\prime \prime }(x) = x^{s - 2} + x^{ - s -
1}.
\end{equation}

Thus we have $\phi _s ^{\prime \prime }(x) > 0$ for all $x > 0$,
and hence, $\phi _s (x)$ is convex for all $x > 0$. Also, we have
$\phi _s (1) = 0$. In view of this we can say that
\textit{J-divergence of type s }is \textit{nonnegative} and
\textit{convex} in the pair of probability distributions $(P,Q)
\in \Gamma _n \times \Gamma _n $.
\end{proof}

\begin{property} \label{prt42} The measure $\mathcal V _s (P\vert \vert Q)$
is monotonically increasing in $s$ for all $s \geqslant
\frac{1}{2}$ and decreasing in $s \leqslant \frac{1}{2}$.
\end{property}

In order to prove the above property, we shall make use the
following lemma.

\begin{lemma} \label{lem41} Let $f:I \subset \mathbb{R}_ + \to
\mathbb{R}$ be a differentiable function and suppose that $f(1) =
f^\prime (1) = 0$, then
\begin{equation}
\label{eq72} f(x)\begin{cases}
 { \geqslant 0,} & {\mbox{if }f\mbox{ is convex}} \\
 { \leqslant 0,} & {\mbox{if }f\mbox{ is concave}} \\
\end{cases}.
\end{equation}
\end{lemma}

\begin{proof} It is well known that if the function $f$ is
convex, then we have the inequality
\begin{equation}
\label{eq73} {f}'(x)(y - x) \leqslant f(y) - f(x) \leqslant
{f}'(y)(y - x),
\end{equation}

\noindent for all $x,y \in \mathbb{R}_ + $. The above inequality
is reversed if $f$ is concave. Take $x = 1$ in the inequality
(\ref{eq73}) and use the fact that $f(1) = f^\prime (1) = 0$ we
get the required result.
\end{proof}

\begin{proof} \textit{of the Property \ref{prt42}}. Let use consider the first order
derivative of the function $\phi _s (x)$ given in (\ref{eq69})
with respect to $s$, we get
\begin{align}
\label{eq74} k_s (x) & = \frac{d}{ds}\left( {\phi _s (x)}
\right)\\
& = \left[ {s(s - 1)} \right]^{ - 2}\left[ {s(s - 1)(x^s - x^{1 -
s})} \right.\ln x\notag\\
& \qquad \left. { + (1 - 2s)(x^s + x^{1 - s} - (x + 1)} \right],
\,\, s \ne 0,1.\notag
\end{align}

Now, calculating the first and second order derivative of the
function $k_s (x)$ with respect to $x$, we get
\begin{align}
k_s ^\prime (x) & = \frac{1}{s^2(1 - s)^2}\left[ {s^2(x^{ - s} -
x^{s - 1}) + (1 - 2s)(x^{ - s} - 1)} \right.\notag\\
& \qquad \left. { + s(s - 1)\left( {sx^{s - 1} + (s - 1)x^{ - s}}
\right)\ln x} \right], \,\, s \ne 0,1\notag
\end{align}

\noindent and
\[
k_s ^{\prime \prime }(x) = (x^{s - 2} - x^{ - s - 1})\ln x.
\]

For all $x > 0$, we can easily check that
\begin{equation}
\label{eq75} {k}''_s (x)\begin{cases}
 { \geqslant 0,} & {s \geqslant \frac{1}{2}} \\
 { \leqslant \mbox{0,}} & {s \leqslant \frac{1}{2}} \\
\end{cases}.
\end{equation}

Since $k_s (1) = {k}'_s (1) = 0$, then using Lemma \ref{lem41}
along with (\ref{eq75}), we have
\begin{equation}
\label{eq76} k_s (x)\begin{cases}
 { \geqslant 0,} & {s \geqslant \frac{1}{2}} \\
 { \leqslant \mbox{0,}} & {s \leqslant \frac{1}{2}} \\
\end{cases}.
\end{equation}

Thus from (\ref{eq76}), we conclude that the function $\phi _s
(x)$ is monotonically increasing in $s$ for all $s \geqslant
\frac{1}{2}$ and monotonically decreasing in $s$ for all $s
\leqslant \frac{1}{2}$. This completes the proof of the property.
\end{proof}

By taking $s = \frac{1}{2}$, 1 and 2, and applying the Property
\ref{prt42}, one gets
\begin{equation}
\label{eq77} h(P\vert \vert Q) \leqslant \frac{1}{8}J(P\vert \vert
Q) \leqslant \frac{1}{16}\Psi (P\vert \vert Q).
\end{equation}

\begin{theorem} \label{the41} The following bounds hold:
\begin{equation}
\label{eq78} \mathcal V _s (P\vert \vert Q) \leqslant E_{\mathcal
V _s } (P\vert \vert Q) \leqslant A_{\mathcal V _s } (r,R),
\end{equation}
\begin{equation}
\label{eq79} \mathcal V _s (P\vert \vert Q) \leqslant B_{\mathcal
V _s } (r,R) \leqslant A_{\mathcal V _s } (r,R),
\end{equation}
\begin{align}
\label{eq80} & \left| {\mathcal V _s (P\vert \vert Q) -
\frac{1}{2}E_{\mathcal V _s } (P\vert \vert Q)} \right|\\
& \leqslant \min \left\{ {\frac{1}{8}\delta _{\mathcal V _s }
(r,R)\chi ^2(P\vert \vert Q),\,\,\frac{1}{12}\left\| {\phi _s
^{\prime \prime \prime }} \right\|_\infty \left| \chi
\right|^3(P\vert \vert Q),\,\,\mathop V\limits_r^R ({\phi
}')V(P\vert \vert Q)} \right\},\notag
\end{align}

\noindent and
\begin{align}
\label{eq81} & \left| {\mathcal V _s (P\vert \vert Q) -
E_{\mathcal V _s
}^\ast (P\vert \vert Q)} \right|\\
& \leqslant \min \left\{ {\frac{1}{8}\delta _{\mathcal V _s }
(r,R)\chi ^2(P\vert \vert Q),\,\,\frac{1}{24}\left\| {\phi _s
^{\prime \prime \prime }} \right\|_\infty \left| \chi
\right|^3(P\vert \vert Q),\,\,\frac{1}{2}\mathop V\limits_r^R
({\phi }')V(P\vert \vert Q)} \right\},\notag
\end{align}

\noindent where
\begin{equation}
\label{eq82} E_{\mathcal V _s } (P\vert \vert Q) =\begin{cases}
 {\sum\limits_{i = 1}^n {(p_i - q_i )\left[ {(s - 1)^{ - 1}\left( {\frac{p_i
}{q_i }} \right)^{s - 1} - s^{ - 1}\left( {\frac{p_i }{q_i }}
\right)^{ -
s}} \right],} } & {s \ne 1} \\\\
 {J(P\vert \vert Q) + \chi ^2(Q\vert \vert P),} & {s = 1} \\
\end{cases},
\end{equation}
\begin{equation}
\label{eq83} E_{\mathcal V _s }^\ast (P\vert \vert Q)
 =\begin{cases}
 {\sum\limits_{i = 1}^n {(p_i - q_i )\left[ {(s - 1)^{ - 1}\left( {\frac{p_i
+ q_i }{2q_i }} \right)^{s - 1} - s^{ - 1}\left( {\frac{p_i + q_i
}{2q_i }}
\right)^{ - s}} \right],} } & {s \ne 0,1} \\\\
 {\Delta (P\vert \vert Q) + 2\, J\left( {\frac{P + Q}{2}\vert \vert Q} \right),}
& {s = 0,1} \\
\end{cases},
\end{equation}
\begin{equation}
\label{eq84} A_{\mathcal V _s } (r,R) = \frac{1}{4}(R -
r)^2\left\{ {L_{s - 2}^{s - 2} (r,R) + L_{ - s - 1}^{ - s - 1}
(r,R)} \right\},
\end{equation}
\begin{equation}
\label{eq85} B_{\mathcal V _s } (r,R) =\begin{cases}
 {[s(s - 1)]^{ - 1}\left[ {\frac{(1 - r)(R^s + R^{1 - s}) + (R - 1)(r^s +
r^{1 - s})}{(R - r)} - 2} \right],} & {s \ne 0,1} \\\\
 {(1 - r)(R - 1)L_{ - 1}^{ - 1} (r,R),} & {s = 0,1} \\
\end{cases},
\end{equation}
\begin{equation}
\label{eq86} \delta _{\mathcal V _s } (r,R) = (R - r)\left[ {(2 -
s)L_{s - 3}^{s - 3} (r,R) + (1 + s)L_{ - s - 2}^{ - s - 2} (r,R)}
\right], \,\,
 - 1 \leqslant s \leqslant 2,
\end{equation}
\begin{equation}
\label{eq87} \left\| {\phi _s ^{\prime \prime \prime }}
\right\|_\infty = (2 - s)r^{s - 3} + (s + 1)r^{ - s - 2}, \,\,
 - 1 \leqslant s \leqslant 2,
\end{equation}

\noindent and
\begin{equation}
\label{eq88} \mathop V\limits_r^R ({\phi }') = \frac{4}{R -
r}A_{\mathcal V _s } (r,R).
\end{equation}
\end{theorem}

\begin{proof} By making some calculations and applying the Theorem \ref{the32}
we get the inequalities (\ref{eq78}) and (\ref{eq79}). Let us
prove now the inequalities (\ref{eq80}) and (\ref{eq81}). The
third order derivative of the function $\phi _s (x)$ is given by
\begin{equation}
\label{eq89} \phi _s ^{\prime \prime \prime }(x) = - \left[ {(2 -
s)x^{s - 3} + (s + 1)x^{ - s - 2}} \right], \,\, x \in (0,\infty
).
\end{equation}

This gives
\begin{equation}
\label{eq90} \phi _s ^{\prime \prime \prime }(x) \leqslant 0, \,\,
\forall \,\, - 1 \leqslant s \leqslant 2.
\end{equation}

From (\ref{eq90}), we can say that the function $\phi ^{\prime
\prime }(x)$ is monotonically decreasing in $x \in (0,\infty )$,
and hence, for all $x \in [r,R]$, we have
\begin{align}
\label{eq91} \delta _{\mathcal V _s } (r,R) & = {\phi }''(r) -
{\phi
}''(R)\\
& = (R - r)\left[ {(2 - s)L_{s - 3}^{s - 3} (r,R) + (1 + s)L_{ - s
- 2}^{ - s - 2} (r,R)} \right], \,\,
 - 1 \leqslant s \leqslant 2.\notag
\end{align}

From (\ref{eq89}), we have
\[
\left| {\phi _s ^{\prime \prime \prime }(x)} \right| = (2 - s)x^{s
- 3} + (s + 1)x^{ - s - 2}, \,\,
 - 1 \leqslant s \leqslant 2.
\]

This gives
\begin{align}
\label{eq92} \left| {{\phi }'''(x)} \right|^\prime & = - \left[
{(s - 2)(s - 3)x^{s - 4} + (s + 1)(s + 2)x^{ - 3 - s}} \right]\\
& \leqslant 0, \,\, - 1 \leqslant s \leqslant 2.\notag
\end{align}

In view of (\ref{eq92}), we can say that the function $\left|
{\phi _s ^{\prime \prime \prime }(x)} \right|$ is monotonically
decreasing in $x \in (0,\infty )$ for $ - 1 \leqslant s \leqslant
2$, and hence, for all $x \in [r,R]$, we have
\begin{equation}
\label{eq93} \left\| {\phi _s ^{\prime \prime \prime }}
\right\|_\infty = \mathop {\sup }\limits_{x \in [r,R]} \left|
{\phi _s ^{\prime \prime \prime }(x)} \right| = (2 - s)r^{s - 3} +
(s + 1)r^{ - s - 2}, \,\,
 - 1 \leqslant s \leqslant 2.
\end{equation}

By applying Theorem \ref{the35} along with the expressions
(\ref{eq91}) and (\ref{eq93}) for the measure (\ref{eq16}) we get
the first two parts of the inequalities (\ref{eq80}) and
(\ref{eq81}). The last part of the inequalities (\ref{eq80}) and
(\ref{eq81}) are obtained by using (\ref{eq88}) and Theorem
\ref{the35}.
\end{proof}

In particular when $s = \frac{1}{2}$, 1 and 2, we get the results
studies in Taneja \cite{tan6, tan7}.

\subsection*{4.2. Bounds on AG and JS -- Divergence of
Type \textit{s}}

Initially we shall give two important properties of \textit{AG and
JS - divergences of type s}.

\begin{property} \label {prt43} The measure $\mathcal W _s (P\vert \vert Q)$
is \textit{nonnegative} and \textit{convex} in the pair of
probability distributions $(P,Q) \in \Gamma _n \times \Gamma _n $
for all $s \in ( - \infty ,\infty )$.
\end{property}

\begin{proof} For all $x > 0$ and $s \in ( - \infty ,\infty
)$, let us consider in (\ref{eq15})
\begin{equation}
\label{eq94} \psi _s (x) =\begin{cases}
 {\left[ {s(s - 1)} \right]^{ - 1}\left[ {\left( {\frac{x^{1 - s} + 1}{2}}
\right)\left( {\frac{x + 1}{2}} \right)^s - \left( {\frac{x +
1}{2}}
\right)} \right],} & {s \ne 0,1} \\\\
 {\frac{x}{2}\ln x - \left( {\frac{x + 1}{2}} \right)\ln \left( {\frac{x +
1}{2}} \right),} & {s = 0} \\\\
 {\left( {\frac{x + 1}{2}} \right)\ln \left( {\frac{x + 1}{2\sqrt x }}
\right),} & {s = 1} \\
\end{cases},
\end{equation}

\noindent then we have $C_f (P\vert \vert Q) = \mathcal W _s
(P\vert \vert Q)$, where $\mathcal W _s (P\vert \vert Q)$ is as
given by (\ref{eq17}).

Moreover,
\begin{equation}
\label{eq95} \psi _s ^\prime (x) =\begin{cases}
 {(s - 1)^{ - 1}\left[ {\frac{1}{s}\left[ {\left( {\frac{x + 1}{2x}}
\right)^s - 1} \right] - \frac{x^{ - s} - 1}{4}\left( {\frac{x +
1}{2}}
\right)^{s - 1}} \right],} & {s \ne 0,1} \\\\
 { - \frac{1}{2}\ln \left( {\frac{x + 1}{2x}} \right),} & {s = 0}
 \\\\
 {1 - x^{ - 1} - \ln x - 2\ln \left( {\frac{2}{x + 1}} \right),} & {s = 1}
\\
\end{cases}
\end{equation}

\noindent and
\begin{equation}
\label{eq96} \psi _s ^{\prime \prime }(x) = \left( {\frac{x^{ - s
- 1} + 1}{8}} \right)\left( {\frac{x + 1}{2}} \right)^{s - 2}.
\end{equation}

Thus we have $\psi _s ^{\prime \prime }(x) > 0$ for all $x > 0$,
and hence, $\psi _s (x)$ is convex for all $x > 0$. Also, we have
$\psi _s (1) = 0$. In view of this we can say that\textit{ AG and
JS -- divergences of type }$s$ is \textit{nonnegative} and
\textit{convex} in the pair of probability distributions $(P,Q)
\in \Gamma _n \times \Gamma _n $.
\end{proof}

\begin{property} \label{prt44} The measure $\mathcal W _s (P\vert \vert Q)$ is
monotonically increasing in $s$ for all $s \geqslant - 1$.
\end{property}

\begin{proof} Let us consider the first order derivative of
(\ref{eq94}) with respect to $s$.
\begin{align}
\label{eq97}
 m_s (x) & = \frac{d}{ds}\left(
{\psi _s (x)} \right)\\
& = - \left[ {s(s - 1)} \right]^{ - 2}\left( {\frac{x + 1}{2}}
\right)^s\left[ {(2s - 1)(x^{1 - s} + x + 2)} \right.\notag \\
& \qquad \qquad \left. { - s(s - 1)(x^{1 - s} + 1)\ln \left(
{\frac{x + 1}{2}} \right)} \right], \,\, s \ne 0,1.\notag
\end{align}

Now, calculating the first and second order derivatives of
(\ref{eq97}) with respect to $x$, we get
\[
m_s ^\prime (x) = \frac{1 - 2s}{s^2(1 - s)^2}\left[ {x^s + x^{1 -
s} - (x + 1)} \right]
 + \frac{1}{s(s - 1)}(x^s - x^{1 - s})\ln x,
\,\, s \ne 0,1
\]

\noindent and
\begin{equation}
\label{eq98} m_s ^{\prime \prime }(x) = \frac{1}{2x^2(x +
1)^2}\left( {\frac{x + 1}{2}} \right)^s\left[ {x^{1 - s}\ln \left(
{\frac{x + 1}{2x}} \right) + x^2\ln \left( {\frac{x + 1}{2}}
\right)} \right],
\end{equation}

\noindent respectively.

Since $(x - 1)^2 \geqslant 0$ for any $x$, this give us
\begin{equation}
\label{eq99} \ln \left( {\frac{x + 1}{2x}} \right) \geqslant \ln
\left( {\frac{2}{x + 1}} \right).
\end{equation}

Now for all $0 < x \leqslant 1$ and for any $s \geqslant - 1$, we
have $x^{1 - s} \geqslant x^2$. This together with (\ref{eq99})
gives
\begin{equation}
\label{eq100} m_s ^{\prime \prime }(x) \geqslant 0,\mbox{ for all
}0 < x \leqslant 1\mbox{ and }s \geqslant - 1.
\end{equation}

Reorganizing (\ref{eq98}), we can write
\[
m_s ^{\prime \prime }(x) = \frac{x^{1 - s}}{2x^2(x + 1)^2}\left(
{\frac{x + 1}{2}} \right)^s\left[ {(x^{1 + s} + 1)\ln \left(
{\frac{x + 1}{2}} \right) - \ln x} \right].
\]

Again for all $x \geqslant 1$ and $s \geqslant - 1$, we have $x^{1
+ s} + 1 \geqslant 2$. This gives
\begin{equation}
\label{eq101} (x^{1 + s} + 1)\ln \left( {\frac{x + 1}{2}} \right)
\geqslant 2\ln \left( {\frac{x + 1}{2}} \right) \geqslant \ln x,
\end{equation}

\noindent where we have used the fact that $(x + 1)^2 \geqslant
4x$ for any $x$.

In view of (\ref{eq101}), we have
\begin{equation}
\label{eq102} m_s ^{\prime \prime }(x) \geqslant 0,\mbox{ for all
}x \geqslant 1\mbox{ and }s \geqslant - 1.
\end{equation}

Combining (\ref{eq100}) and (\ref{eq102}), we have
\begin{equation}
\label{eq103} m_s ^{\prime \prime }(x) \geqslant 0,\mbox{ for all
}x > 0\mbox{ and }s \geqslant - 1.
\end{equation}

Since $m_s (1) = m_s ^\prime (1) = 0$, then (\ref{eq103}) together
with Lemma \ref{lem41} complete the required proof.
\end{proof}

By taking $s = - 1$, 0, $\frac{1}{2}$, 1 and 2, and applying
Property \ref{prt44}, one gets
\begin{equation}
\label{eq104} \frac{1}{4}\Delta (P\vert \vert Q) \leqslant
I(P\vert \vert Q) \leqslant 4\,d(P\vert \vert Q) \leqslant
T(P\vert \vert Q) \leqslant \frac{1}{16}\Psi (P\vert \vert Q).
\end{equation}

\begin{theorem} \label{the42} The following bounds hold:
\begin{equation}
\label{eq105} 0 \leqslant \mathcal W _s (P\vert \vert Q) \leqslant
E_{\mathcal W _s } (P\vert \vert Q) \leqslant A_{\mathcal W _s }
(r,R),
\end{equation}
\begin{equation}
\label{eq106} 0 \leqslant \mathcal W _s (P\vert \vert Q) \leqslant
B_{\mathcal W _s } (r,R) \leqslant A_{\mathcal W _s } (r,R),
\end{equation}
\begin{align}
\label{eq107} & \left| {\mathcal W _s (P\vert \vert Q) -
\frac{1}{2}E_{\mathcal W _s } (P\vert \vert Q)} \right|\\
& \leqslant \min \left\{ {\frac{1}{8}\delta _{\mathcal W _s }
(r,R)\chi ^2(P\vert \vert Q),\,\,\frac{1}{12}\left\| {\psi _s
^{\prime \prime \prime }} \right\|_\infty \left| \chi
\right|^3(P\vert \vert Q),\,\,\mathop V\limits_r^R ({\psi
}')V(P\vert \vert Q)} \right\},\notag
\end{align}

\noindent and
\begin{align}
\label{eq108} & \left| {\mathcal W _s (P\vert \vert Q) -
E_{\mathcal W _s
}^\ast (P\vert \vert Q)} \right|\\
& \leqslant \min \left\{ {\frac{1}{8}\delta _{\mathcal W _s }
(r,R)\chi ^2(P\vert \vert Q),\,\,\frac{1}{24}\left\| {\psi _s
^{\prime \prime \prime }} \right\|_\infty \left| \chi
\right|^3(P\vert \vert Q),\,\,\frac{1}{2}\mathop V\limits_r^R
({\psi }')V(P\vert \vert Q)} \right\},\notag
\end{align}

\noindent where
\begin{align}
\label{eq109} & E_{\mathcal W _s } (P\vert \vert Q)\\
& =\begin{cases}
 {\frac{1}{2}\sum\limits_{i = 1}^n {(p_i - q_i )\left\{ {(s - 1)^{ -
1}\left( {\frac{p_i^{1 - s} + q_i^{1 - s} }{2}} \right)\left(
{\frac{p_i + q_i }{2}} \right)^{s - 1}\left. { - s^{ - 1}\left(
{\frac{p_i + q_i }{2p_i
}} \right)^s} \right\} ,} \right.} } & {s \ne 0,1} \\\\
 {J\left( {\frac{P + Q}{2}\vert \vert P} \right),} & {s = 0} \\\\
 {\frac{1}{4}\left[ {\chi ^2(Q\vert \vert P) - J(P\vert \vert Q)} \right] +
J\left( {\frac{P + Q}{2}\vert \vert Q} \right),} & {s = 1} \\
\end{cases},\notag
\end{align}
\begin{align}
\label{eq110} & E_{\mathcal W _s }^\ast (P\vert \vert Q)\\
& =\begin{cases}
 {\left( {\frac{1}{2}} \right)^{s + 1}\sum\limits_{i = 1}^n {(p_i - q_i
)\left\{ {(s - 1)^{ - 1}\left[ {\left( {\frac{p_i + 3q_i }{p_i +
q_i }}
\right)^{s - 1}} \right.} \right.} } & \\
 \qquad \qquad {\left. {\left. { + \left( {\frac{p_i + 3q_i }{2q_i }} \right)^{s -
1}} \right] - s^{ - 1}\left( {\frac{p_i + 3q_i }{2q_i }}
\right)^s} \right\}
,} & {s \ne 0,1} \\\\
 {2J\left( {\frac{P + Q}{2}\vert \vert \frac{P + 3Q}{4}} \right),} & {s = 0}
\\\\
 {\frac{1}{4}\Delta (P\vert \vert Q) - \frac{1}{2}J\left( {\frac{P +
Q}{2}\vert \vert Q} \right) + 2J\left( {\frac{P + 3Q}{4}\vert
\vert Q}
\right),} & {s = 1} \\
\end{cases},\notag
\end{align}
\begin{align}
\label{eq111} A_{\mathcal W _s } (r,R) & = \frac{(R -
r)^2}{16}\left[ {\frac{1}{rR}L_{s - 1}^{s - 1} \left( {\frac{r +
1}{2r},\frac{R +
1}{2R}} \right)} \right.\\
& \qquad \left. { - \frac{1}{2rR}L_{s - 2}^{s - 2} \left( {\frac{r
+ 1}{2r},\frac{R + 1}{2r}} \right) + \frac{1}{2}L_{s - 2}^{s - 2}
\left( {\frac{r + 1}{2},\frac{R + 1}{2}} \right)} \right],\notag
\end{align}

\begin{align}
\label{eq112} & B_{\mathcal V _s } (r,R)\\
 & = \begin{cases}
 {[s(s - 1)]^{ - 1}\left\{ {\frac{1}{R - r}\left[ {(1 - r)\left( {\frac{R^{1
- s} + 1}{2}} \right)\left( {\frac{R + 1}{2}} \right)^s} \right.}
\right.} &
\\
 {\left. {\left. {\qquad \qquad + (R - 1)\left( {\frac{r^{1 - s} + 1}{2}}
\right)\left( {\frac{r + 1}{2}} \right)^s} \right] - 1} \right\}
,} & {s \ne
0,1} \\\\
 {\frac{1}{2(R - r)}\left\{ {(1 - r)\left[ {R\ln R - (1 + R)\ln \left(
{\frac{R + 1}{2}} \right)} \right]} \right.} & \\
 {\left. {\qquad + (R - 1)\left[ {r\ln r - (r + 1)\ln \left( {\frac{r +
1}{2}} \right)} \right]} \right\} ,} & {s = 0} \\\\
 {\frac{1}{2}\left\{ {(1 - rR)L_{ - 1}^{ - 1} (r + 1,R + 1) + \ln \left[
{\frac{(r + 1)(R + 1)}{8}} \right]} \right\},} & {s = 1} \\
\end{cases},\notag
\end{align}
\begin{align}
\label{eq113} \delta _{\mathcal W _s } (r,R) & = \left( {\frac{r^{
- s -
1} + 1}{8}} \right)\left( {\frac{r + 1}{2}} \right)^{s - 2} \\
& \qquad - \left( {\frac{R^{ - s - 1} + 1}{8}} \right)\left(
{\frac{R + 1}{2}} \right)^{s - 2}, \,\,
 - 1 \leqslant s \leqslant 2,\notag
\end{align}
\begin{align}
\label{eq114} \left\| {\psi _s ^{\prime \prime \prime }}
\right\|_\infty & = \frac{1}{2(r^3 + 1)}\left( {\frac{r + 1}{2}}
\right)^s \times\\
& \qquad \times \left[ {3r^{ - s - 1} + (s + 1)r^{ - s - 2} + (2 -
s)} \right], \,\,
 - 1 \leqslant s \leqslant 2,\notag
\end{align}

\noindent and
\begin{equation}
\label{eq115} \mathop V\limits_r^R ({\psi }') = \frac{4}{R -
r}A_{\mathcal W_s} (r,R).
\end{equation}
\end{theorem}

\begin{proof} By making some calculations and applying
Theorem \ref{the32} we get the inequalities (\ref{eq105}) and
(\ref{eq106}). Now, we shall prove the inequalities (\ref{eq107})
and (\ref{eq108}). The third order derivative of the function
$\psi _s (x)$ is given by
\begin{align}
\label{eq116} \psi _s ^{\prime \prime \prime }(x) & = -
\frac{1}{2(x + 1)^3}\left( {\frac{x + 1}{2}} \right)^s\times\\
& \qquad \qquad \times \left[ {3x^{ - s - 1} + (s + 1)x^{ - 2 - s}
+ (2 - s)} \right], \,\, x \in (0,\infty ).\notag
\end{align}

This gives
\begin{equation}
\label{eq117} \psi _s ^{\prime \prime \prime }(x) \leqslant 0,
\,\, - 1 \leqslant s \leqslant 2.
\end{equation}

From (\ref{eq117}), we can say that the function $\phi ^{\prime
\prime }(x)$ is monotonically decreasing in $x \in (0,\infty )$,
and hence, for all $x \in [r,R]$, we have
\begin{align}
\label{eq118} \delta _{^{\mathcal W _s }} (r,R) & = {\psi }''(r) -
{\psi }''(R)\\
& = \left( {\frac{r^{ - s - 1} + 1}{8}} \right)\left( {\frac{r +
1}{2}} \right)^{s - 2} \notag\\
& \qquad - \left( {\frac{R^{ - s - 1} + 1}{8}} \right)\left(
{\frac{R + 1}{2}} \right)^{s - 2}, \,\,
 - 1 \leqslant s \leqslant 2.\notag
\end{align}

Again, from (\ref{eq116}), we have
\begin{align}
\left| {\psi _s ^{\prime \prime \prime }(x)} \right| & =
\frac{1}{2(x + 1)^3}\left( {\frac{x + 1}{2}}
\right)^s \times\\
 & \qquad \times \left[ {3x^{ - s - 1} + (s
+ 1)x^{ - 2 - s} + (2 - s)} \right], \,\, x \in (0,\infty ), \,\,
 - 1 \leqslant s \leqslant 2. \notag
\end{align}

This gives
\begin{align}
\label{eq119} \left| {\psi _s ^{\prime \prime \prime }(x)} \right|
& = - \frac{x^{1 - s}}{2(x + 1)^4} \times\\
& \qquad \times \left\{ {x^{1 - s}\left[ {12x^2 + 8(s + 1)x + (s +
1)(s + 2)} \right] + x^4(s - 2)(s - 3)} \right\}\notag\\
& \leqslant 0, \,\, - 1 \leqslant s \leqslant 2. \notag
\end{align}

In view of (\ref{eq119}), we can say that the function $\left|
{\psi _s ^{\prime \prime \prime }} \right|$ is monotonically
decreasing in $x \in (0,\infty )$ for $ - 1 \leqslant s \leqslant
2$, and hence, for all $x \in [r,R]$, we have
\begin{align}
\label{eq120} \left\| {\psi _s ^{\prime \prime \prime }}
\right\|_\infty & = \mathop {\sup }\limits_{x \in [r,R]} \left|
{\psi _s ^{\prime \prime \prime }(x)} \right|\\
& = \frac{1}{2(r + 1)^3}\left( {\frac{r + 1}{2}} \right)^s\left[
{3r^{ - s - 1} + (s + 1)r^{ - 2 - s} + (2 - s)} \right], \,\,
 - 1 \leqslant s \leqslant 2.\notag
\end{align}

By applying the Theorem \ref{the35} along with the expressions
(\ref{eq118}) and (\ref{eq120}) for the measure (\ref{eq17}) we
get the first two parts of the bounds (\ref{eq107}) and
(\ref{eq108}). The last part of the bounds (\ref{eq107}) and
(\ref{eq108}) follows in view of (\ref{eq115}) and Theorem
\ref{the35}.
\end{proof}

In particular when $s = - 1$, 0, 1 and 2, we get the results
studied in Taneja \cite{tan6, tan7}.

\section{Relations Among Generalized Relative
Divergence Measures}

In this section, we shall apply the Theorem \ref{the36} to obtain
inequalities among the measures (\ref{eq16}) and (\ref{eq17}).\\

Let us consider
\begin{equation}
\label{eq121} g_{(\psi _s ,\phi _t )} (x) = \frac{{\psi }''_s
(x)}{\phi _t ^{\prime \prime }(x)} = \frac{x^{ - 1 - s} + 1}{8(x^{
- t - 1} + x^{t - 2})}\left( {\frac{x + 1}{2}} \right)^{s - 2},
\,\, x \in (0,\infty )
\end{equation}

\noindent where ${\psi }''_s (x)$ and ${\phi }''_t (x)$ are as
given by (\ref{eq96}) and (\ref{eq71}) respectively.

From (\ref{eq121}) one has
\begin{align}
\label{eq122} {g}'_{(\psi _s ,\phi _t )} (x) & = \frac{1}{8x(x +
1)(x^{ - t - 1} + x^{t - 2})^2}\left( {\frac{x + 1}{2}} \right)^{s
- 2} \times\\
& \qquad \times \left\{ {x^{ - t - 1}} \right.\left[ {x^{ - s -
1}\left( {(t - 2)x + (t - s)} \right) + (t + s - 1)x + (t + 1)}
\right]\notag\\
& \qquad \qquad \left. { - x^{t - 2}\left[ {x^{ - s - 1}\left( {(t
+ 1)x + (t +
s - 1)} \right) + (t - s)x + (t - 2)} \right]} \right\}\notag\\
& = \frac{1}{8x(x + 1)(x^{ - t - 1} + x^{t - 2})^2}\left( {\frac{x
+ 1}{2}} \right)^{s - 2}\left[ {(t - s)(x^{ - s - t - 2} - x^{t -
1})} \right.\notag \\
& \qquad  + (t - 2)(x^{ - s - t - 1} - x^{t - 2})\notag\\
& \qquad \qquad \left. { + (t + s - 1)(x^{ - t} - x^{t - s - 3}) +
(t + 1)(x^{ - t - 1} - x^{t - s - 2})} \right].\notag\\
& = \frac{1}{8x(x + 1)(x^{ - t - 1} + x^{t - 2})^2}\left( {\frac{x
+ 1}{2}} \right)^{s - 2}\left\{ {(t - s)\left[ {x^{t - 1}(x^{ - 2t
- s - 1} - 1)} \right]} \right.\notag\\
& \qquad  + (t - 2)\left[ {x^{t - 2}(x^{ - s - 2t + 1} - 1)}
\right] + (t + s - 1)\left[ {x^{t - s - 3}(x^{s - 2t + 3} - 1)}
\right]\notag\\
& \qquad \qquad \left. { + (t + 1)\left[ {x^{t - s - 2}(x^{s - 2t
+ 1} - 1)} \right]} \right\}.\notag
\end{align}

From the above expression we observe that it is difficult to know
the nature of the expression (\ref{eq122}) with respect to the
parameters $s$ and $t$. Here below we shall study the above
expression for some particular cases of the parameters $s$ and
$t$.

\bigskip
\subsection*{5.1. Inequalities Among $\mathcal W _s (P\vert \vert
Q)$ and $J(P\vert \vert Q)$} Take $t = 1$ in (\ref{eq121}) and
(\ref{eq122}), we get
\begin{equation}
\label{eq123} g_{(\psi _s ,\phi _1 )} (x) = \frac{{\psi }''_s
(x)}{\phi _1 ^{\prime \prime }(x)} = \frac{x^{ - 1 - s} + 1}{8(x^{
- 2} + x^{ - 1})}\left( {\frac{x + 1}{2}} \right)^{s - 2}, \,\, x
\in (0,\infty )
\end{equation}

\noindent and
\begin{equation}
\label{eq124} {g}'_{(\psi _s ,\phi _1 )} (x) = - \frac{x}{8(x +
1)^2}\left( {\frac{x + 1}{2}} \right)^{s - 2}\left[ {2(x^{ - s} -
1) + (s - 1)x(x^{ - s - 2} - 1)} \right] \,\, .
\end{equation}

\noindent respectively.

Using the fact that
\begin{equation}
\label{eq125} x^k\begin{cases}
 { \geqslant 1,} & {x \geqslant 1,\,\,k > 0\mbox{ or }x \leqslant
1,\,\,k < 0} \\
 { \leqslant 1,} & {x \leqslant 1,\,\,k > 0\mbox{ or }x \geqslant
1,\,\,k < 0} \\
\end{cases}.
\end{equation}

\noindent we can write the expression (\ref{eq124}) as follow:
\begin{equation}
\label{eq126} {g}'_{(\psi _s ,\phi _1 )} (x)\begin{cases}
 { \leqslant 0,} & {(x \geqslant 1,\,\, - 2 \leqslant s \leqslant
0)\mbox{ or }(x \leqslant 1,\,\,s \geqslant 1)} \\
 { \geqslant 0,} & {(x \leqslant 1,\,\, - 2 \leqslant s \leqslant
0)\mbox{ or }(x \geqslant 1,\,\,s \geqslant 1)} \\
\end{cases}.
\end{equation}

From (\ref{eq126}), we conclude that the function $g_{(\psi _s
,\phi _1 )} (x)$ is monotonically decreasing (resp. increasing) in
$(1,\infty )$ and increasing (resp. decreasing) in $(0,1)$ for all
$ - 2 \leqslant s \leqslant 0$ (resp. $s \geqslant 1)$. This gives
\begin{equation}
\label{eq127} M = \mathop {\sup }\limits_{x \in (0,\infty )}
g_{(\psi _s ,\phi _1 )} (x) = g_{(\psi _s ,\phi _1 )} (1) =
\frac{1}{8}, \,\,
 - 2 \leqslant s \leqslant 0,
\end{equation}

\noindent and
\begin{equation}
\label{eq128} m = \mathop {\inf }\limits_{x \in (0,\infty )}
g_{(\psi _s ,\phi _1 )} (x) = g_{(\psi _s ,\phi _1 )} (1) =
\frac{1}{8}, \,\, s \geqslant 1.
\end{equation}

By the application of the inequality (\ref{eq60}) given in Theorem
\ref{the36} with the expressions (\ref{eq127}) and (\ref{eq128}),
we conclude the following inequality among the measures $\mathcal
W _s (P\vert \vert Q)$ and $J(P\vert \vert Q)$:
\begin{equation}
\label{eq129} \mathcal W _s (P\vert \vert Q)\begin{cases}
 { \leqslant \frac{1}{8}J(P\vert \vert Q),} & { - 2 \leqslant s \leqslant 0}
\\\\
 { \geqslant \frac{1}{8}J(P\vert \vert Q),} & {s \geqslant 1} \\
\end{cases}.
\end{equation}

In particular the expression (\ref{eq129}) lead us to to the
inequality
\begin{equation}
\label{eq130} I(P\vert \vert Q) \leqslant \frac{1}{8}J(P\vert
\vert Q) \leqslant T(P\vert \vert Q).
\end{equation}

\bigskip
\subsection*{5.2. Inequalities Among $\mathcal W _s (P\vert \vert
Q)$ and  $h(P\vert \vert Q)$} Take $t = 1/2$ \, in (\ref{eq121})
and (\ref{eq122}), we get
\begin{equation}
\label{eq131} g_{(\psi _s ,\phi _{1 / 2} )} (x) = \frac{{\psi
}''_s (x)}{\phi _{1 / 2} ^{\prime \prime }(x)} = \frac{x^{ - 1 -
s} + 1}{16x^{ - 3 / 2}}\left( {\frac{x + 1}{2}} \right)^{s - 2},
\,\, x \in (0,\infty )
\end{equation}

\noindent and
\begin{align}
\label{eq132} {g}'_{(\psi _s ,\phi _{1 / 2} )} (x) & = -
\frac{\sqrt x }{32(x + 1)^2}\left( {\frac{x + 1}{2}} \right)^{s -
2}\times\\
& \qquad \qquad \times \left[ {x^{ - s - 1}\left( {3x + (2s - 1)}
\right)
- (2s - 1)x - 3} \right]\notag\\
& = - \frac{\sqrt x }{32(x + 1)^2}\left( {\frac{x + 1}{2}}
\right)^{s - 2}\left[ {3(x^{ - s} - 1) + (2s - 1)x(x^{ - s - 2} -
1)} \right].\notag
\end{align}

\noindent respectively.

Again in view of (\ref{eq125}) we can write
\begin{equation}
\label{eq133} {g}'_{(\psi _s ,\phi _{1 / 2} )} (x)\begin{cases}
 { \leqslant 0,} & {(x \geqslant 1,\,\, - 2 \leqslant s \leqslant
0)\mbox{ or }(x \leqslant 1,\,\,s \geqslant \frac{1}{2})} \\\\
 { \geqslant 0,} & {(x \leqslant 1,\,\, - 2 \leqslant s \leqslant
0)\mbox{ or }(x \geqslant 1,\,\,s \geqslant \frac{1}{2})} \\
\end{cases}.
\end{equation}

From (\ref{eq133}), we conclude that the function $g_{(\psi _s
,\phi _{1 / 2} )} (x)$ is monotonically decreasing (resp.
increasing) in $(1,\infty )$ and increasing (resp. decreasing) in
$(0,1)$ for all $ - 2 \leqslant s \leqslant 0$ (resp. $s \geqslant
\frac{1}{2})$. This gives
\begin{equation}
\label{eq134} M = \mathop {\sup }\limits_{x \in (0,\infty )}
g_{(\psi _s ,\phi _{1 / 2} )} (x) = g_{(\psi _s ,\phi _{1 / 2} )}
(1) = \frac{1}{8}, \,\,
 - 2 \leqslant s \leqslant 0,
\end{equation}

\noindent and
\begin{equation}
\label{eq135} m = \mathop {\inf }\limits_{x \in (0,\infty )}
g_{(\psi _s ,\phi _{1 / 2} )} (x) = g_{(\psi _s ,\phi _{1 / 2} )}
(1) = \frac{1}{8}, \,\, s \geqslant \frac{1}{2}.
\end{equation}

By the application of the inequality (\ref{eq60}) given in Theorem
\ref{the36} with the expressions (\ref{eq134}) and (\ref{eq135}),
we conclude the following inequality among the measures $\mathcal
W _s (P\vert \vert Q)$ and $h(P\vert \vert Q)$:
\begin{equation}
\label{eq136} \mathcal W _s (P\vert \vert Q)\begin{cases}
 { \leqslant h(P\vert \vert Q),} & { - 2 \leqslant s \leqslant 0}
 \\\\
 { \geqslant h(P\vert \vert Q),} & {s \geqslant \frac{1}{2}} \\
\end{cases}.
\end{equation}

In particular this gives
\begin{equation}
\label{eq137} I(P\vert \vert Q) \leqslant h(P\vert \vert Q)
\leqslant T(P\vert \vert Q).
\end{equation}

The expressions (\ref{eq77}), (\ref{eq104}), (\ref{eq130}) and
(\ref{eq137}) together give an inequality among the six
\textit{symmetric divergence measures} as follows:
\begin{equation}
\label{eq138} \frac{1}{4}\Delta (P\vert \vert Q) \leqslant
I(P\vert \vert Q) \leqslant h(P\vert \vert Q) \leqslant
\frac{1}{8}J(P\vert \vert Q) \leqslant T(P\vert \vert Q) \leqslant
\frac{1}{16}\Psi (P\vert \vert Q).
\end{equation}

\begin{remark} $(i)$ In view of (\ref{eq138}) and
(\ref{eq10}), we have the following inequality:
\begin{equation}
\label{eq139} \frac{1}{4}\Delta (P\vert \vert Q) \leqslant
I(P\vert \vert Q) \leqslant h(P\vert \vert Q) \leqslant
\frac{1}{8}J(P\vert \vert Q) \leqslant T(P\vert \vert Q) \leqslant
\frac{1}{4}J(P\vert \vert Q).
\end{equation}

$(ii)$ It is well known \cite{lec} that
\begin{equation}
\label{eq140} \frac{1}{4}\Delta (P\vert \vert Q) \leqslant
h(P\vert \vert Q) \leqslant \frac{1}{2}\Delta (P\vert \vert Q).
\end{equation}

In view of (\ref{eq138}) and (\ref{eq140}), we have the following
inequality:
\begin{equation}
\label{eq141} \frac{1}{4}\Delta (P\vert \vert Q) \leqslant
I(P\vert \vert Q) \leqslant h(P\vert \vert Q) \leqslant
\frac{1}{2}\Delta (P\vert \vert Q).
\end{equation}
\end{remark}

\bigskip
\subsection*{5.3. Inequalities Among $\mathcal W _s (P\vert \vert
Q)$, $\mathcal V _t (P\vert \vert Q)$ and $\Psi (P\vert \vert Q)$}
For $t = 2$, we have
\[
\mathcal V _2 (P\vert \vert Q) = \frac{1}{2}\Psi (P\vert \vert Q)
\]

\noindent and
\[
\mathcal W _2 (P\vert \vert Q) = \frac{1}{16}\Psi (P\vert \vert
Q).
\]

Since the measure $\mathcal W _s (P\vert \vert Q)$ is
monotonically increasing in $s$ for all $s \geqslant - 1$. This
gives
\begin{equation}
\label{eq142} \mathcal W _s (P\vert \vert Q)\begin{cases}
 { \leqslant \frac{1}{16}\Psi (P\vert \vert Q),} & { - 1 \leqslant s
\leqslant 2} \\\\
 { \geqslant \frac{1}{16}\Psi (P\vert \vert Q),} & {s \geqslant 2} \\
\end{cases}.
\end{equation}

Also, in view of monotonicity of $\mathcal V _t (P\vert \vert Q)$
with respect to $t$, we have
\begin{equation}
\label{eq143} \mathcal V _t (P\vert \vert Q)\begin{cases}
 { \leqslant \frac{1}{2}\Psi (P\vert \vert Q),} & {\frac{1}{2} \leqslant t
\leqslant 2} \\\\
 { \geqslant \frac{1}{2}\Psi (P\vert \vert Q),} & {s \leqslant \frac{1}{2}}
\\
\end{cases}.
\end{equation}

\bigskip
\subsection*{5.4. Inequalities Among $\mathcal V _t (P\vert \vert
Q)$ and $\Delta (P\vert \vert Q)$} Take $s = -1$ in (\ref{eq121})
and (\ref{eq122}), we get
\begin{equation}
\label{eq144} g_{(\psi _{ - 1} ,\phi _t )} (x) = \frac{{\psi }''_{
- 1} (x)}{\phi _t ^{\prime \prime }(x)} = \frac{1}{2(x^{ - t - 1}
+ x^{t - 2})(x + 1)}, \,\, x \in (0,\infty )
\end{equation}

\noindent and
\begin{align}
\label{eq145} {g}'_{(\psi _{ - 1} ,\phi _t )} (x) & = \frac{2}{x(x
+ 1)^4(x^{ - t - 1} + x^{t - 2})^2}\times\\
& \qquad \times\left[ {(t - 2)(x^{ - t} - x^{t - 2}) + (t + 1)(x^{
- t - 1} - x^{t - 1})} \right]\notag\\
& = \frac{2}{x(x + 1)^4(x^{ - t - 1} + x^{t -
2})^2}\times\notag\\
& \qquad \left[ {(t - 2)x^{t - 2}(x^{2(1 - t)} - 1) + (t + 1)x^{t
- 1}(x^{ - 2t} - 1)} \right].\notag
\end{align}

\noindent respectively.

Again in view of (\ref{eq125}) we can write
\begin{equation}
\label{eq146} {g}'_{(\psi _{ - 1} ,\phi _t )} (x)\begin{cases}
 { \leqslant 0,} & {x \geqslant 1,\mbox{ }t \geqslant 0 \,\, \mbox{or} \,\, t \leqslant -1} \\
 { \geqslant 0,} & {x \leqslant 1,\mbox{ }t \geqslant 0 \,\, \mbox{or} \,\, t \leqslant -1} \\
\end{cases}.
\end{equation}

From (\ref{eq146}) we conclude that the function $g_{(\psi _{ - 1}
,\phi _t )} (x)$ is monotonically decreasing (resp. increasing) in
$(1,\infty )$ and increasing (resp. decreasing) in $(0,1)$ for all
$t \geqslant 0 \,\, \mbox{or} \,\, t \leqslant -1$. This gives
\begin{equation}
\label{eq147} M = \mathop {\sup }\limits_{x \in (0,\infty )}
g_{(\psi _{ - 1} ,\phi _t )} (x) = g_{(\psi _{ - 1} ,\phi _t )}
(1) = \frac{1}{8}, \,\, t \geqslant 0 \,\, \mbox{or} \,\, t
\leqslant -1.
\end{equation}

By the application of the inequality (\ref{eq60}) given in Theorem
\ref{the36} with the expression (\ref{eq147}), we conclude the
following inequality among the measures $\mathcal V _t (P\vert
\vert Q)$ and $\Delta (P\vert \vert Q)$:
\begin{equation}
\label{eq148} \Delta (P\vert \vert Q) \leqslant
\frac{1}{2}\mathcal V _t (P\vert \vert Q), \,\, t \geqslant 0 \,\,
\mbox{or} \,\, t \leqslant -1.
\end{equation}

\bigskip
\subsection*{5.5. Inequalities Among $\mathcal V _t (P\vert \vert
Q)$ and $I(P\vert \vert Q)$} Take $s = 0$ in (\ref{eq121}) and
(\ref{eq122}), we get
\begin{equation}
\label{eq149} g_{(\psi _0 ,\phi _t )} (x) = \frac{{\psi }''_0
(x)}{\phi _t ^{\prime \prime }(x)} = \frac{1}{4x(x + 1)(x^{ - t -
1} + x^{t - 2})}, \,\, x \in (0,\infty )
\end{equation}

\noindent and
\begin{align}
\label{eq150} {g}'_{(\psi _0 ,\phi _t )} (x) & = \frac{1}{2x^2(x +
1)^2(x^{ - t - 1} + x^{t - 2})^2}\times\\
& \qquad \times \left[ {(t - 1)(x^{ - t} - x^{t - 2}) + t(x^{ - t
- 1} - x^{t - 1})} \right]\notag\\
& = \frac{2}{2x^2(x + 1)^2(x^{ - t - 1} + x^{t -
2})^2}\times\notag\\
& \qquad \left[ {(t - 1)x^{t - 2}(x^{2(1 - t)} - 1) + tx^{t -
1}(x^{ - 2t} - 1)} \right].\notag
\end{align}

\noindent respectively.

Again in view of (\ref{eq125}) we can write
\begin{equation}
\label{eq151} {g}'_{(\psi _0 ,\phi _t )} (x)\begin{cases}
 { \leqslant 0,} & {x \geqslant 1,\mbox{ }} \\
 { \geqslant 0,} & {x \leqslant 1,\mbox{ }} \\
\end{cases}.
\end{equation}

\noindent for all $t \in ( - \infty ,\infty )$.

From (\ref{eq151}) we conclude that the function $g_{(\psi _0
,\phi _t )} (x)$ is monotonically decreasing in $(1,\infty )$ and
increasing in $(0,1)$ for all $ - \infty < t < \infty $. This
gives
\begin{equation}
\label{eq152} M = \mathop {\sup }\limits_{x \in (0,\infty )}
g_{(\psi _0 ,\phi _t )} (x) = g_{(\psi _0 ,\phi _t )} (1) =
\frac{1}{8}, \,\,
 - \infty < t < \infty .
\end{equation}

By the application of the inequality (\ref{eq60}) with the
expression (\ref{eq152}), we conclude the following inequality
among the measures $\mathcal V _s (P\vert \vert Q)$ and $I(P\vert
\vert Q)$:
\begin{equation}
\label{eq153} I(P\vert \vert Q) \leqslant \frac{1}{8}\mathcal V _t
(P\vert \vert Q), \,\,
 - \infty < t < \infty .
\end{equation}

\bigskip
\subsection*{5.6. Inequalities Among $ \mathcal V _t (P\vert \vert
Q)$ and  $T(P\vert \vert Q)$} Take $s = 1$ in (\ref{eq121}) and
(\ref{eq122}), we get
\begin{equation}
\label{eq154} g_{(\psi _1 ,\phi _t )} (x)  = \frac{{\psi }''_1
(x)}{\phi _t ^{\prime \prime }(x)} = \frac{x^2 + 1}{4x^2(x +
1)(x^{ - t - 1} + x^{t - 2})}, \,\, x \in (0,\infty )
\end{equation}

\noindent and
\begin{align}
\label{eq155} {g}'_{(\psi _1 ,\phi _t )} (x) & = \frac{1}{4x^4(x +
1)^2(x^{ - t - 1} + x^{t - 2})^2}\times\\
& \qquad \times\left[ {t(x^{2 - t} - x^{t - 2}) + (t + 1)(x^{ - t
+ 1} - x^{t - 1})} \right.\notag\\
&\qquad \qquad \left. { + (t - 2)(x^{ - t} - x^t) + (t - 1)(x^{ -
t - 1}
- x^{t + 1})} \right]\notag\\
& = \frac{1}{4x^4(x + 1)^2(x^{ - t - 1} + x^{t -
2})^2}\times\notag\\
& \qquad \times\left[ {tx^{t - 2}(x^{2(2 - t)} - 1) + (t + 1)x^{t
- 1}(x^{2(1 - t)} - 1)} \right.\notag\\
& \qquad \qquad \left. { + (t - 2)x^t(x^{ - 2t} - 1) + (t - 1)x^{t
+ 1}(x^{ - 2(t + 1)} - 1)} \right].\notag
\end{align}

\noindent respectively.

Again in view of (\ref{eq125}) we can write
\begin{equation}
\label{eq156} {g}'_{(\psi _1 ,\phi _t )} (x)\begin{cases}
 { \geqslant 0,} & {(x \geqslant 1,\mbox{ }0 \leqslant t \leqslant 1),\mbox{
}(x \leqslant 1,\mbox{ }t \geqslant \mbox{2} \,\, \mbox{or} \,\, t\leqslant -1)} \\
 { \leqslant 0,} & {(x \leqslant 1,\mbox{ }0 \leqslant t \leqslant 1),\mbox{
}(x \geqslant 1,\mbox{ }t \geqslant \mbox{2} \,\, \mbox{or} \,\, t\leqslant -1)} \\
\end{cases}.
\end{equation}

From (\ref{eq156}) we conclude that the function $g_{(\psi _1
,\phi _t )} (x)$ is monotonically decreasing (resp. increasing) in
$(1,\infty )$ and increasing (resp. decreasing) in $(0,1)$ for all
$0 \leqslant t \leqslant 1$(resp. $t \geqslant 2$ or $t \leqslant
- 1)$. This gives
\begin{equation}
\label{eq157} m = \mathop {\sup }\limits_{x \in (0,\infty )}
g_{(\psi _1 ,\phi _t )} (x) = g_{(\psi _1 ,\phi _t )} (1) =
\frac{1}{8}, \,\, 0 \leqslant t \leqslant 1,
\end{equation}

\noindent and
\begin{equation}
\label{eq158} M = \mathop {\inf }\limits_{x \in (0,\infty )}
g_{(\psi _1 ,\phi _t )} (x) = g_{(\psi _1 ,\phi _t )} (1) =
\frac{1}{8}, \,\, t \geqslant 2\mbox{ or }t \leqslant - 1.
\end{equation}

By the application of the inequality (\ref{eq60}) with the
expressions (\ref{eq157}) and (\ref{eq158}), we conclude the
following inequality among the measures $\mathcal V _t (P\vert
\vert Q)$ and $T(P\vert \vert Q)$:
\begin{equation}
\label{eq159} T(P\vert \vert Q)\begin{cases}
 { \leqslant \frac{1}{8}\mathcal V _t (P\vert \vert Q),} & {t \geqslant 2\mbox{
or }t \geqslant - 1} \\\\
 { \geqslant \frac{1}{8}\mathcal V _t (P\vert \vert Q),} & {0 \leqslant t
\leqslant 1} \\
\end{cases}.
\end{equation}

\bigskip
\subsection*{5.7. Inequalities Among $\mathcal W _s (P\vert \vert
Q)$ and $\mathcal V _s (P\vert \vert Q)$} When $t = s$, we shall
obtain the results in two different ways, one by using the Theorem
\ref{the36}, and another applying Jensen's inequality.

(i) Take $t = s$ in (\ref{eq121}) and (\ref{eq122}), we get
\begin{equation}
\label{eq160} g_{(\psi _s ,\phi _s )} (x) = \frac{{\psi }''_s
(x)}{{\phi }''_s (x)} = \frac{x^{ - 1 - s} + 1}{8(x^{ - s - 1} +
x^{s - 2})}\left( {\frac{x + 1}{2}} \right)^{s - 2}, \,\, x \in
(0,\infty )
\end{equation}

\noindent and
\begin{align}
\label{eq161} {g}'_{(\psi _s ,\phi _t )} (x) & = \frac{1}{8x(x +
1)(x^{ - s - 1} + x^{s - 2})^2}\left( {\frac{x + 1}{2}} \right)^{s
- 2}\times\\
& \qquad \times \left[ {(s + 1)x^{ - 2}(x^{1 - s} - 1) + (s -
2)x^{s - 2}(x^{1 - 3s} - 1)} \right.\notag\\
& \qquad \qquad \left. { + (2s - 1)x^{ - 3}(x^{3 - s} - 1)}
\right]\notag
\end{align}

\noindent respectively.

Again in view of (\ref{eq125}) we can write
\begin{equation}
\label{eq162} {g}'_{(\psi _s ,\phi _s )} (x)\begin{cases}
 { \geqslant 0,} & {(x \geqslant 1,\mbox{ }\frac{1}{2} \leqslant s \leqslant
1),\mbox{ }(x \leqslant 1,\mbox{ }s \geqslant 2\mbox{ or }s
\leqslant - 1)}
\\
 { \leqslant 0,} & {(x \leqslant 1,\mbox{ }\frac{1}{2} \leqslant s \leqslant
1),\mbox{ }(x \geqslant 1,\mbox{ }s \geqslant 2\mbox{ or }s
\leqslant - 1)}
\\
\end{cases}.
\end{equation}

From (\ref{eq162}) we conclude that the function $g_{(\psi _s
,\phi _s )} (x)$ is monotonically increasing (resp. decreasing) in
$(1,\infty )$ and decreasing (resp. increasing) in $(0,1)$ for all
$\frac{1}{2} \leqslant s \leqslant 1$(resp. $s \geqslant 2$ or $s
\leqslant - 1)$. This gives
\begin{equation}
\label{eq163} m = \mathop {\inf }\limits_{x \in (0,\infty )}
g_{(\psi _s ,\phi _s )} (x) = g_{(\psi _s ,\phi _s )} (1) =
\frac{1}{8}, \,\, \frac{1}{2} \leqslant s \leqslant 1,
\end{equation}

\noindent and
\begin{equation}
\label{eq164} M = \mathop {\sup }\limits_{x \in (0,\infty )}
g_{(\psi _s ,\phi _s )} (x) = g_{(\psi _s ,\phi _s )} (1) =
\frac{1}{8}, \,\, s \geqslant 2\mbox{ or }s \leqslant - 1.
\end{equation}

By the application of the inequality (\ref{eq60}) with the
expressions (\ref{eq163}) and (\ref{eq164}), we conclude the
following inequality among the measures $\mathcal W _s (P\vert
\vert Q)$ and $\mathcal V _s (P\vert \vert Q)$:
\begin{equation}
\label{eq165} \mathcal W _s (P\vert \vert Q)\begin{cases}
 { \leqslant \frac{1}{8}\mathcal V _s (P\vert \vert Q),} & {s \geqslant 2\mbox{
or }s \geqslant - 1} \\\\
 { \geqslant \frac{1}{8}\mathcal V _s (P\vert \vert Q),} & {\frac{1}{2} \leqslant
s \leqslant 1} \\
\end{cases}.
\end{equation}

Some particular cases of (\ref{eq165}) can be seen in (518) or
(\ref{eq139}).\\

(ii) By applying Jensen's inequality we can easily check that
\begin{equation}
\label{eq166} \frac{p_i^s + q_i^s }{2}\begin{cases}
 { \leqslant \left( {\frac{p_i + q_i }{2}} \right)^s,} & {0 < s < 1}
 \\\\
 { \geqslant \left( {\frac{p_i + q_i }{2}} \right)^s,} & {s > 1\mbox{ or }s
< 0} \\
\end{cases},
\end{equation}

\noindent for all $i = 1,2,...,n$, where $P,Q \in \Gamma _n $.

Multiplying (\ref{eq166}) by $p_i^{1 - s} $, summing over all $i =
1,2,...,n$ and simplifying, we get
\begin{equation}
\label{eq167} \frac{1}{2}\left( {1 + \sum\limits_{i = 1}^n {p_i^{1
- s} q_i^s } } \right)\begin{cases}
 { \leqslant \sum\limits_{i = 1}^n {p_i^{1 - s} \left( {\frac{p_i + q_i
}{2}} \right)^s} ,} & {0 < s < 1} \\\\
 { \geqslant \sum\limits_{i = 1}^n {p_i^{1 - s} \left( {\frac{p_i + q_i
}{2}} \right)^s} ,} & {s > 1\mbox{ or s} < 0} \\
\end{cases}.
\end{equation}

Similarly, we can write
\begin{equation}
\label{eq168} \frac{1}{2}\left( {1 + \sum\limits_{i = 1}^n {q_i^{1
- s} p_i^s } } \right)\begin{cases}
 { \leqslant \sum\limits_{i = 1}^n {q_i^{1 - s} \left( {\frac{p_i + q_i
}{2}} \right)^s} ,} & {0 < s < 1} \\\\
 { \geqslant \sum\limits_{i = 1}^n {q_i^{1 - s} \left( {\frac{p_i + q_i
}{2}} \right)^s} ,} & {s > 1\mbox{ or }s < 0} \\
\end{cases}.
\end{equation}

Adding (\ref{eq167}) and (\ref{eq168}) and making some
adjustments, we get
\begin{align}
\label{eq169} & \frac{1}{4}\left( {\sum\limits_{i = 1}^n {p_i^{1 -
s} q_i^s +
q_i^{1 - s} p_i^s - 2} } \right) \\
& \qquad \begin{cases}
 { \leqslant \sum\limits_{i = 1}^n {\left( {\frac{p_i^{1 - s} + q_i^{1 - s}
}{2}} \right)\left( {\frac{p_i + q_i }{2}} \right)^s - 1} ,} & {0
< s < 1}
\\
 { \geqslant \sum\limits_{i = 1}^n {\left( {\frac{p_i^{1 - s} + q_i^{1 - s}
}{2}} \right)\left( {\frac{p_i + q_i }{2}} \right)^s - 1} ,} & {s
> 1\mbox{
or }s < 0} \\
\end{cases}.\notag
\end{align}

Since, $s(s - 1) < 0$ for all $0 < s < 1$ and $s(s - 1) > 0$ for
all $s > 1$ or $s < 0$. This together with (\ref{eq170}) proves
that
\begin{equation}
\label{eq170} \mathcal V _s (P\vert \vert Q) \geqslant 4\mathcal W
_s (P\vert \vert Q), \,\,
 - \infty < s < \infty ,
\end{equation}

\noindent where for $s = 0$ and $s = 1$, the result is obtained by
the continuity of the measures with respect to the parameter $s$.

In view of (\ref{eq165}) and (\ref{eq170}), we can write
\begin{equation}
\label{eq171} \frac{1}{8}\mathcal V _s (P\vert \vert Q) \leqslant
\mathcal W _s (P\vert \vert Q) \leqslant \frac{1}{4} \mathcal V _s
(P\vert \vert Q), \,\, \frac{1}{2} \leqslant s \leqslant 1.
\end{equation}

\begin{remark}
\begin{itemize}
\item[(i)] For $s = 1$ in (\ref{eq171}), we get a part of the
inequality (\ref{eq139}). For $s = \frac{1}{2}$ in (\ref{eq171}),
we get an interesting bound on a new measure given by
(\ref{eq181}) in terms of \textit{Hellinger's discrimination}:
\begin{equation}
\label{eq172} \frac{1}{4}h(P\vert \vert Q) \leqslant d(P\vert
\vert Q) \leqslant \frac{1}{2}h(P\vert \vert Q).
\end{equation}

\item[(ii)] In particular for $s = \frac{1}{2}$ and $t = 1$ in
(\ref{eq121}), we can easily show that
\begin{equation}
\label{eq182} 4\,d(P\vert \vert Q) \leqslant \frac{1}{8}J(P\vert
\vert Q).
\end{equation}
The inequalities (\ref{eq138}) together with (\ref{eq172}) and
(\ref{eq182}) give the following inequalities among some
particular cases of the measures (\ref{eq16}) and (\ref{eq17}):
\begin{align}
\label{eq183} & \frac{1}{4}\Delta (P\vert \vert Q) \leqslant
I(P\vert \vert Q) \leqslant h(P\vert \vert Q) \leqslant
4\,d(P\vert
\vert Q)\\
&\qquad \leqslant \frac{1}{8}J(P\vert \vert Q) \leqslant T(P\vert
\vert Q) \leqslant \frac{1}{16}\Psi (P\vert \vert Q).\notag
\end{align}
\end{itemize}
\end{remark}


\begin{thebibliography}{99}
\bibitem{pea} K. PEARSON, On the Criterion that a given system of
eviations from the probable in the case of correlated system of
variables is such that it can be reasonable supposed to have
arisen from random sampling, \textit{Phil. Mag.},
\textbf{50}(1900), 157-172.
\end{thebibliography}
\end{document}